\documentclass[oneside,dvipdfmx]{amsart}
\usepackage{amsmath, amssymb,amsthm,cite}
\usepackage{graphicx}
\usepackage{color}
\newtheorem{theorem}{Theorem}[section]
\newtheorem{proposition}[theorem]{Proposition}
\newtheorem{lemma}[theorem]{Lemma}
\newtheorem{corollary}[theorem]{Corollary}

\newtheorem{question}{Question}
\newtheorem{rem}{Remark}
\newtheorem{example}{Example}

\theoremstyle{definition}

\newtheorem{main}{Theorem}

\def\nbd{neighborhood } 
 
\def\T{\mathbb{T} } 
\def\F{\mathcal{F} } 
\def\G{\mathcal{G} }

\def\Fh{\hat{\mathcal{F}} } 
 
\def\Q{\mathbb{Q} } 
\def\R{\mathbb{R} } 
\def\Z{\mathbb{Z} } 
\makeatletter
\@namedef{subjclassname@2020}{%
  \textup{2020} Mathematics Subject Classification}
\makeatother

\title[On low-codimensional foliations]{On low-codimensional foliations}

\author{Tomoo YOKOYAMA}
\date{\today}
\address{Applied Mathematics and Physics Division, Gifu University, Yanagido 1-1, Gifu, 501-1193, Japan\\}
\email{tomoo@gifu-u.ac.jp}
\subjclass[2020]{57R30 (primary), 57M60,54B15,37C10,37E30 (secondary).}
\keywords{Foliations, group actions, decompositions, flows, homeomorphisms}

\thanks{This work was partially supported by JSPS Kakenhi Grant Number 20K03583}

\begin{document}
\maketitle

\begin{abstract}
Decompositions on manifolds appear in various geometric structures. Necessary and sufficient conditions for quotient spaces of decompositions to be manifolds are widely characterized. We characterize necessary and sufficient conditions to be $k$-manifolds $(k = 1, 2)$, which generalize characterizations in the codimension-$k$ cases for the leaf spaces of foliations, the orbit spaces of group-actions, decomposition spaces of upper semi-continuous decompositions, and leaf class spaces of Riemannian foliation with regular closure. To prove such characterizations, we generalize a characterization of upper semi-continuity for decomposition into one for a class decomposition. In addition, we completely characterize upper semi-continuity for class decompositions of homeomorphisms on orientable compact surfaces and of homeomorphisms isotopic to identity on non-orientable compact surfaces. Furthermore, a flow on a connected compact $3$-manifold whose class decomposition is upper semi-continuous is ``almost $k$ dimensional'' $(k=0, 1, 2, 3)$ or has ``complicated'' minimal sets. 
\end{abstract}

\section{Introduction}

Decompositions on manifolds are appeared in various geometric structures.  
For instance, necessary and sufficient conditions for quotient spaces of decompositions to be manifolds are widely characterized (e.g. \cite{D2,D3,DW,DW2,EMS,E,E2,V,V2}). 
We review properties of quotient spaces of decompositions from points of decomposition theory, of foliation theory, of group actions, and of dynamical systems and generalize some results as follows.

\subsection{Decompositions}

In \cite{D2}, for any integer $n \geq 1$, it is shown that the decomposition space $M/\mathcal{F}$ of an upper semi-continuous decomposition $\mathcal{F}$ on an $n$-manifold $M$ whose elements are connected closed $(n-1)$-manifolds is a $1$-manifold, where the decomposition space $M/\mathcal{F}$ is a quotient space $M/\sim_{\F}$ defined by $x \sim_{\F} y$ if $\mathcal{F}(x) = \mathcal{F}(y)$. 
Here $\mathcal{F}(x)$ is the element of $\mathcal{F}$ containing $x$.  
Moreover, for any integer $n \geq 2$, the decomposition space $M/\mathcal{F}$ of an upper semi-continuous decomposition $\mathcal{F}$ on an $n$-manifold $M$ whose elements are connected closed $(n-2)$-manifolds is a surface (i.e. $2$-manifold) \cite{D3,DW,DW2}. 
In addition, the decomposition space of an upper semi-continuous decomposition on an $n$-manifold whose elements are connected closed $(n-3)$-manifolds is $3$-dimensional if the elements are orientable \cite{D4}. 
On the other hand, the decomposition space of an upper semi-continuous decomposition on an $n$-manifold whose elements are connected closed $(n-3)$-manifolds need be neither a manifold nor even a generalized manifold (for instance, see  \cite[Example 3.3]{DW3}). 
Here a generalized $n$-manifold is a finite-dimensional ANR having the local homology of $\R^n$. 
However, for any integers $k \geq 2$ and $n \geq 2k$, the decomposition space $M/\mathcal{F}$ of an upper semi-continuous decomposition on an $n$-manifold $M$ whose elements are connected closed orientable $n-k$-manifolds whose $l$th integral \v{C}ech homology groups are trivial for any $l = 1, 2, \ldots , k-1$ is a generalized $k$-manifold unless the covering dimension of $M/\mathcal{F}$ is infinite (see \cite[Theorem~3.1]{DW3} for detail).
The decomposition space of an upper semi-continuous decomposition on a $(2m+1)$-manifold ($m \in \mathbb{Z}_{>0}$) whose elements are conpacta having the shape of homology $m$-spheres is a generalized $m+1$-manifold if the covering dimension of $M/\mathcal{F}$ is finite and the discontinuity set does not locally separate the decomposition space (see \cite[Theorem~1]{snyder1991} for detail). 

Recall that the class space of a decomposition is the $T_0$-tification of the decomposition space. 
In other words, the class space $X/\hat{\mathcal{F}}$ of a decomposition $\mathcal{F}$ on a topological space $X$ is a quotient space $X/\sim_{\hat{\F}}$ defined by $x \sim_{\hat{\F}} y$ if $\overline{\mathcal{F}(x)} = \overline{\mathcal{F}(y)}$, where $\overline{A}$ is the closure of a subset $A$.  
Recall that a decomposition is pointwise almost periodic if and only if the set of all closures of elements of $\mathcal{F}$ also is a decomposition. 
The element class $\hat{L}$ of an element $L$ of a decomposition $\F$ is defined by $\hat{L} := \{ y \in M \in \overline{\F(y)} = \overline{L} \}$. 
The set of element classes is a decomposition, called the class decomposition and denoted by $\hat{\mathcal{F}}$. 
Note that the class space  of a decomposition on a topological space  corresponds to the decomposition space of the class decomposition. 
A decomposition $\F$ on a topological space $M$ is $R$-closed if the subset $R = \{ (x, y) \in M \times M \mid y \in \overline{\F(x)} \}$ is a closed subset. 
%

We partially generalize a characterization of upper semi-continuity for the decomposition space of a decompositions into compact and connected elements of a locally compact Hausdorff space \cite[Remark after Theorem~4.1]{E2}, into one for a class decomposition into arbitrary elements of a (compact) $T_4$ space as follows. 
 
\begin{main}\label{th:usc}
The following conditions are equivalent for a decomposition $\mathcal{F}$ of a $T_4$ space $X$: 
\\
$(1)$ The decomposition $\mathcal{F}$ is $R$-closed.  \\
$(2)$ The class space $X/\hat{\mathcal{F}}$ is Hausdorff. 
$(3)$ The decomposition $\mathcal{F}$ is of characteristic $0$. 
\\
Moreover, each of the following conditions implies the conditions $(1)$-$(3)$, and if $X$ is compact then each of the following conditions is equivalent to one of the conditions $(1)$-$(3)$: 
\\
$(4)$ The class decomposition $\hat{\mathcal{F}}$ is upper semi-continuous.  
\\
$(5)$ The decomposition $\mathcal{F}$ is $L$-stable $(\mathrm{i.e.}$ For any neighborhood $U$ of each element $\hat{L} \in \hat{\mathcal{F}}$, there is an open ${\F}$-invariant neighborhood of $\hat{L}$ contained in $U)$. 
\\
$(6)$ The decomposition $\mathcal{F}$ is weakly almost periodic {\rm(i.e.} The quotient map $p \colon  X \to X/\hat{\mathcal{F}}$ is closed{\rm)}. 

In any case, the decomposition $\mathcal{F}$ is pointwise almost periodic. 
\end{main}

The compactness of the space $X$ in the previous theorem is necessary, because 
any $R$-closed decomposition $\F = \{ X \}$ for a non-compact $T_4$ space $X$ is weakly almost periodic but $\hat{\mathcal{F}} = \F$ is not upper semi-continuous. 

\subsection{Foliations}

Since higher codimension foliations are more flexible than codimension one foliations, 
several studies of higher codimension foliations are restricted to some classes (e.g. a class of compact foliations \cite{EMS,E,E2,ER,Mi,V,V4,V2,Z}, a class of compact Hausdorff foliations \cite{CH,dN,F,FM,H,NR,V5}). 
Recall that a compact foliation $\F$ on a manifold $M$  is a foliation with all leaves compact and a compact Hausdorff foliation is a compact foliation whose leaf space $M/\F$ is a Hausdorff topological space with the quotient topology, where the leaf space $M/\F$ is a quotient space collapsing leaves into singletons (i.e. a quotient space $M/\sim_{\F}$ defined by $x \sim_{\F} y$ if the leaf $\mathcal{F}(x)$ containing $x$ corresponds to the leaf $\mathcal{F}(y)$ containing $y$ (i.e. $\mathcal{F}(x) = \mathcal{F}(y)$)). 
It is known that the leaf space of a continuous codimension two compact foliation of a compact manifold is a compact orbifold \cite{EMS,E,E2,V,V2}. 
In particular, such a foliation is compact Hausdorff. 
On the other hand, there are non-Hausdorff compact foliations and non-Hausdorff flows each of whose orbits is compact for codimension $q > 2$ \cite{EV,S,V3}. 
Recall that the leaf class $\hat{L}$ of a leaf $L$ of a foliation $\F$ is defined by $\hat{L} := \{ y \in M \in \overline{\F(y)} = \overline{L} \}$. 
Hausdorff separation axiom of a manifold $M$ implies that the leaf space $M/\F$ of a compact foliation $\F$ of $M$ corresponds to the leaf class space $M/\hat{\mathcal{F}}$, where the leaf class space $M/\hat{\mathcal{F}}$ is a quotient space collapsing leaf classes into singletons (i.e. quotient space $M/\sim_{\hat{\F}}$ defined by $x \sim_{\hat{\F}} y$ if $\hat{\mathcal{F}}(x) = \hat{\mathcal{F}}(y)$ (i.e. $\overline{\F(x)} = \overline{\F(y)}$)). 
Note that the leaf class space $M/{\hat{\F}}$ is the $T_0$-tification of the leaf space $M/\mathcal{F}$. 
It is shown \cite{Y} that a foliation of a compact manifold is compact Hausdorff if and only if it is compact and $R$-closed, and that the leaf class space of a foliation of a compact manifold is Hausdorff if and only if the foliation is $R$-closed. 
This means that $R$-closed foliations are generalizations of compact Hausdorff foliations. 
%

For a non-negative integer $k$, a foliation $\F$ is codimension-$k$-like if all but finitely many leaf classes of $\F$ are codimension $k$ submanifolds without boundaries and the finite exceptions are  leaf classes which are closed subsets and whose codimensions are more than $k$.   
Note that each Riemannian foliation with regular closure on a compact manifold is codimension-$k$-like and $R$-closed \cite{molino1988riemannian}. 
Using the characterization of $R$-closedness of decompositions, we generalize the fact that the leaf (class) space of a continuous codimension one/two compact foliation of a compact manifold is a compact orbifold as follows.  

\begin{main}\label{main:a}
Let $\F$ be a foliation of a connected compact manifold $M$ and $\hat{\mathcal{F}}$ the set of leaf classes of $\F$. 
Then the following properties hold:
\\
$(1)$ If the foliation $\F$ is codimension-one-like, then the class space $M/\hat{\mathcal{F}}$ is either a closed interval or a circle if and only if $\mathcal{F}$ is $R$-closed.  
\\
$(2)$ If the foliation $\F$ is codimension-two-like, then the class space $M/\hat{\mathcal{F}}$ is a compact surface with corners if and only if $\mathcal{F}$ is $R$-closed. 
\end{main}

The previous theorem is a generalization of necessary and sufficient conditions for the orbit spaces of codimension-$k$ foliations to be $k$-manifolds $(k = 1, 2)$, and is a generalization of sufficient conditions to be orbifolds for Riemannian foliations with regular closure on compact manifolds whose leaf closures are codimensional $k$ $(k = 1, 2)$\cite[Proposition~5.2]{molino1988riemannian}. 
In fact, Theorem~\ref{main:a} implies that, for any codimension $k$ foliation $\F$ of a compact manifold $M$ ($k =1,2$), the leaf space $M/\mathcal{F}$ is a compact $k$-manifold with or without boundary if and only if the foliation $\mathcal{F}$ is (compact) Hausdorff. 



\subsection{Group actions}

A group-action $G$ is $R$-closed if the set $\mathcal{F} _G$ of orbits is an $R$-closed decomposition. 
In \cite{AGW,Ha}, it is shown that the following properties are equivalent for each action of a finitely generated group $G$ on either a compact zero-dimensional space or a graph $X$: (1) Pointwise recurrent;  (2) Pointwise almost periodic; (3) $R$-closed. 
In addition, the following are equivalent for countable group $G$ on a local dendrite:  (1) Pointwise almost periodic;  (2) Hausdorff for the orbit class space; (3) $R$-closed \cite{Marzougui2019min}. 
Furthermore, the equivalence holds for equicontinuously generated groups on locally compact zero-dimensional spaces \cite{Re}. 
On the other hand, there are group-actions which are not $R$-closed but pointwise almost periodic, and there are group-actions which are not pointwise almost periodic but pointwise recurrent \cite{AGW}. 
As pointed in \cite{salem2018group}, in general, the equivalence of equicontinuity and $R$-closedness for group-actions as mentioned in \cite[Exercise~6. p.46]{auslander1988minimal}, does not hold even for distal homeomorphisms (i.e. continuous $\Z$-action) (see Introduction of the paper \cite{salem2018group} for details). 
%

For a non-negative integer $k$, a group-action $G$ is codimension-$k$-like if all but finitely many orbit closures of $\F_G$ are codimension $k$ connected submanifolds without boundaries and the finite exceptions are connected closed subsets each of whose codimension is more than $k$.  
For a property P, a group is virtually P if there is a finite index subgroup has the property P. 
Similarly, we call that a group-action $G$ is a virtually codimension-$k$-like group-action by a finite index normal subgroup if there is a finite index normal subgroup whose action is codimension-$k$-like. 
We have the following characterization of $R$-closedness for group-actions. 

\begin{main}\label{main:b}
Let $G$ be a virtually codimension-$k$-like group-action by a finite index normal subgroup on a compact connected manifold. 
Then the following properties holds:
\\
$(1)$ When $k = 1$, the orbit class space is either a closed interval or a circle if and only if the group-action is $R$-closed. 
\\
$(2)$ When $k = 2$, the orbit class space is a surface with corners if and only if the group-action is $R$-closed. 
\end{main}

The previous theorem is a generalization of necessary and sufficient conditions for the orbit spaces of codimension-$k$ group-actions to be $k$-manifolds $(k = 1, 2)$. 
In fact, Theorem~\ref{main:b} implies that, for any codimension $k$ group-action of a compact manifold $M$ ($k =1,2$), the orbit space $M/G$ is a compact $k$-manifold with or without boundary if and only if the group-action is Hausdorff.

In \cite{ES}, they have shown that if a continuous mapping $f$ of a topological space $X$ in itself is pointwise recurrent (resp. pointwise almost periodic) then $f^k$ for each positive integer $k$ is pointwise recurrent (resp. pointwise almost periodic). 
This result also holds for group-action cases (see  \cite[Theorem~2.24, 4.04, and 7.04]{GH}). 
Analogously, the following result holds for $R$-closed group-action, which generalize \cite[Lemma~1.1]{Y3} to group-action cases. 

\begin{main}\label{main:c}
Let $G$ be a group-action on a locally connected compact Hausdorff space $X$ and $H$ a finite index normal subgroup. 
Then $G$ is $R$-closed if and only if so is $H$. 
\end{main}


\subsection{Dynamical systems}
The following properties are equivalent for a flow and a homeomorphism on a compact metric space \cite{El,Y}: (1) $R$-closed; (2) $D$-stable; (3) Lyapunov stable; (4) Weakly almost periodic in the sense of Gottschalk; (5) Characteristic 0. 
Note weakly almost periodicity in the sense of Gottschalk implies pointwise almost periodicity by definition. 
In \cite{KT,PX}, minimal sets of surface homeomorphisms are classified.  
In \cite{GKT,KT,KT2018}, topologies of invariant sets of area-preserving surface homeomorphism homotopic to the identity are described. 
In these papers, a filling of an invariant set is one of the key tools to describe topologies of invariant sets. 
Using such filling constructions, we also study codimension two dynamical systems (e.g. surface homeomorphisms, flows in $3$-manifolds). 
\subsubsection{Complete characterizations of $R$-closedness for homeomorphisms}
We completely characterize $R$-closedness for homeomorphisms on orientable compact surfaces in \S~\ref{sec:ori} and for homeomorphisms isotopic to the identity on non-orientable compact surfaces in \S~\ref{sec:nonori}. 
For instance, characterizations for homeomorphisms on tori and orientable hyperbolic surfaces are stated as follows. 

\begin{main}\label{main:d} 
A toral homeomorphism $f \colon  \T^2 \to \T^2$ is $R$-closed if and only if it satisfies one of the following statements:
\\
$(1)$ The homeomorphism $f$ is minimal.
\\
$(2)$ The homeomorphism $f$ is periodic $(\mathrm{i.e.}$ there is a positive integer $n$ such that $f^n$ is identical$)$.
\\
$(3)$ The homeomorphism $f$ is pointwise almost periodic and each minimal set is a finite disjoint union of essential circloids. 
\end{main}

\begin{main}\label{main:e} 
A homeomorphism $f \colon S \to S$ of an orientable compact surface S whose genus is more than one is $R$-closed if and only if it is periodic.
\end{main}

\subsubsection{Pentachotomy of $R$-closed flows on connected compact $3$-manifolds}
For an $R$-closed flow $v$ on a manifold $M$, define the quotient space $M/\sim_{\hat{v}}$ as follows: $x \sim_{\hat{v}} y$ if their orbit closures correspond to each other (i.e. $\overline{O(x)} = \overline{O(y)}$), called the (orbit) class space and denote by $M/\hat{v}$. 
Then 
we state the following pentachotomy that an $R$-closed flow on a connected compact $3$-manifold is either ``almost three-dimensional'',  ``almost two-dimensional'', ``almost one-dimensional'', ``almost zero dimensional'', or with ``complicated'' minimal sets.

\begin{main}\label{main:f} 
For a flow $v$ on a connected compact $3$-manifold $M$, if $v$ is $R$-closed, then one of the following statements holds: 
\\
$(1)$ The flow $v$ is identical $(\mathrm{i.e.}$ $M/ \hat{v}$ is the whole $3$-manifold $M$ $)$.   
\\
$(2)$ The flow $v$ is minimal $(\mathrm{i.e.}$ $M/ \hat{v}$ is a singleton $)$. 
\\
$(3)$ The orbit class space $M/ \hat{v}$ is either $[0,1]$ or $\mathbb{S}^1$ and each interior point of $M/ \hat{v}$ is two-dimensional.   
\\
$(4)$ The periodic point set $\mathrm{Per}(v)$ is open dense, the quotient $\mathrm{Per}(v)/v$ is an open connected surface, the end set of $\mathrm{Per}(v)/v$ is zero-dimensional and homeomorphic to $\mathrm{Sing}(v)$, and $M = \mathrm{Sing}(v) \sqcup \mathrm{Per}(v)$, 
where $\sqcup$ denotes a disjoint union.  
\\
$(5)$ There is a two-dimensional minimal set which is not a suspension of a circloid.  

Conversely, if one of conditions $(1)$--$(4)$ holds, then $v$ is $R$-closed. 
\end{main}

Using several examples, we explain the details of the previous theorem as follows:  
Non-trivial $R$-closed flows on a connected compact $3$-manifolds consisting of closed orbits are of type $(4)$. 
The suspension flows of diffeomorphisms which are irrational rotations with respect to an axis on $\mathbb{S}^2$ or $\mathbb{T}^2$ are of type $(3)$. 
Moreover there are suspension flows of homeomorphisms of type $(3)$ with minimal sets which are not circles but circloids. 
In fact, recall that there is the Remage's $R$-closed homeomorphism in \cite[Example 1]{R} (resp. Herman's $C^{\infty}$ diffeomorphism \cite{He}) on $\mathbb{S}^2$. 
Moreover each of them  can be modified into a toral $R$-closed homeomorphism with pseudo-circles (resp. $C^{\infty}$ diffeomorphism with non-locally-connected continua) as minimal sets \cite{Y2}. 
Using $R$-closed diffeomorphisms on annuli with minimal sets which are circloids but neither circles nor pseudo-circles \cite{N,W}, one can construct $R$-closed diffeomorphisms with such a minimal set on $\mathbb{S}^2$ and $\mathbb{T}^2$. 
Then the suspension flows of such diffeomorphisms are also of type $(3)$. 

Recall that a minimal set is trivial if it is either a closed orbit or the whole manifold. 
The author would like to know the following statement holds or not: A flow on a connected compact $3$-manifold is $R$-closed if and only if one of conditions $(1)$--$(4)$ in Theorem~\ref{main:f} holds.  
In other words, one consider the following question.
\begin{question}\label{q:01}
Are there $R$-closed flows on a connected compact $3$-manifolds with nontrivial minimal sets which are not suspensions of circloids?
\end{question}


\subsection{Contents}
The present paper consists of nine sections. 
In the next section, as preliminaries, we introduce fundamental notions coming from general topology, foliation theory, group-actions, and dynamical systems. 
Moreover, examples of codimension-one-like and codimension-two-like $R$-closed flows are included in the end of the section. 
Section 3 describes the key tool Theorem~\ref{th:usc} that is a characterization of $R$-closedness for  decompositions, which is a generalization of one of upper semi-continuous decomposition. 
Section 4 states results about the $R$-closedness which are analogous results for recurrence and pointwise almost periodicity. 
In particular, Theorem~\ref{main:c} is proved. 
In Section 5, a detailed analysis of cases leads to Theorem~\ref{main:a} and Theorem~\ref{main:b}. 
Section 6 develops the properties of pointwise almost periodicity and $R$-closedness. 
In section 7, completely characterize $R$-closedness for homeomorphisms on orientable compact surfaces. 
In particular, Theorem~\ref{main:d} and Theorem~\ref{main:e} are proved. 
Section 8 states a complete characterizations of $R$-closedness for homeomorphisms isotopic to identity on non-orientable compact surfaces. 
In the final section, Theorem~\ref{main:f} is proved. It states the following five possibilities that an $R$-closed flow on a connected compact $3$-manifold either is ``almost $k$ dimensional'' $(k=0, 1, 2, 3)$, or has ``complicated'' minimal sets. 

\section{Preliminaries}

\subsection{General topology}
\subsubsection{Fundamental concepts}
A singleton is a set consisting of one point. 
A topological space is $T_0$ if for any points $x \neq y \in X$, there is an open subset $U$ of $X$ such that $|\{x, y \} \cap U| =1$.  
The $T_0$-tification (or Kolmogorov quotient) of a topological space is a quotient space $X/\sim$ defined by $x \sim y$ if $\overline{\{ x \}} = \overline{\{ y \}}$. 
Notice that a topological space is $T_0$ if and only if its $T_0$-tification corresponds to it. 
A topological space is $T_1$ if each singleton of it is closed. 
Recall that a topological space is regular if any pair of a closed subset $C$ and a point $p \notin C$ has disjoint open neighborhoods, and a topological space is $T_3$ if it is Hausdorff and regular. 
A topological space is normal if any pair of two disjoint closed subsets has disjoint open neighborhoods, and that a topological space is $T_4$ if it is $T_1$ and normal.
Notice that a topological space is $T_4$ if and only if it is Hausdorff and normal, and that a compact space is $T_4$ if and only if it is Hausdorff. 

By a continuum, we mean a nonempty compact connected metrizable space which is not a singleton. 
%
By a Cantor set, we mean a nonempty totally disconnected perfect compact metrizable space. 
Note that any Cantor set is homeomorphic to the Cantor ternary set. 
A surface is a $2$-manifold with or without boundary and a surface with corners is a $2$-manifold locally modeled by $[0, 1]^2$.
A surface with boundary corresponds to a surface with corners as a topological manifold. 
A $k$-orbifold is a space locally modeled on the quotients of $k$-dimensional Euclidean space $\R^k$ by finite groups. 
Notice that the underlying spaces of $2$-orbifolds (resp. $1$-orbifolds) are $2$-manifolds with corners (resp. $1$-manifolds with or without boundary) (cf. \cite[a description after Theorem~2.3, p.24]{CHK}). 

\subsubsection{Concepts on surfaces}
By an open (resp. closed) disk, we mean that a topological space which is homeomorphic to a unit open (resp. closed) disk in $\R^2$. 
Note that an open disk is homeomorphic to any nonempty simply connected open proper subset in a plane or a sphere. 
Let $D$ be an open disk. 
A subset $B \subset D$ is bounded in $D$ if there is a homeomorphism $f \colon  \R^2 \to D$ such that $f^{-1}(B)$ is bounded in $\R^2$. 
Note that a subset $B \subset D$ is bounded in $D$ if and only if there is a closed disk in $D$ containing $B$. 
Define a filling $\mathrm{Fill}_{D}(A) \subseteq D$ of a subset $A \subsetneq D$ as follows:  
$p \in \mathrm{Fill}_{D}(A)$ if either $p \in A$ or there is a bounded open disk in $D$ whose boundary is contained in $A$ and which contains $p$. 
For an open subset $A$ of $D$, 
we have that $ p \in \mathrm{Fill}_{D}(A)$ if and only if 
there is a simple closed curve $\gamma \subseteq A$ which bounds an open disk containing $p$. 
Note that we do not assume that a subset $A$ is connected and that the concepts of fillings are used in \cite{GKT,JK,KT,KT2018,PX} but they are slightly different. 

A subset of a surface $S$ is essential if it is not null homotopic in $S^*$, where $S^*$ is the resulting closed surface from $S$ by collapsing all boundary components into singletons.
Here a boundary component is a connected component of the boundary. 
Note that a one-dimensional cell complex $\gamma$ in a compact surface $S$ is essential if and only if $\gamma$ is not null homologous with respect to the relative homology group $H_1(S, \partial S; \Z)$, where $\partial S$ is the boundary of the surface $S$ as the boundary of a manifold. 

\subsubsection{Circloids}
A continuum $C$ of a topological space $X$ is annular if it has a neighborhood $U \subset X$ homeomorphic to an open annulus $\mathbb{S}^1 \times (0,1)$ such that $U - C$ has exactly two components each of which is homeomorphic to an open annulus, where $A - B$ is used instead of the set difference $A \setminus B$ when $A \subseteq B$. 
We call such $U$ an annular neighborhood of $C$.
We say a subset $C \subset X$ is a circloid if it is an annular continuum and does not contain any strictly smaller annular continuum as a subset.

\subsubsection{Dimensions and codimensions}
By dimension, we mean Lebesgue covering dimension.  
By Urysohn's theorem, the Lebesgue covering dimension, the large inductive dimension, and the small inductive dimension correspond in separable metrizable spaces (cf. \cite[Theorem~IV. 1]{Nag}). 
It is known that the dimension of a finite union of closed subsets is the maximal of the dimension of such subsets (cf. \cite[Corollary~1.2.6]{Co}). 
For a number $k \in \mathbb{Z}_{\geq 0}$, a subset $A$ of a topological space $X$ is codimension $k$ if the difference between the dimensions of $X$ and $A$ is $k$. 

\subsection{Decompositions}
\subsubsection{Decompositions and their decomposition spaces}
By a decomposition, we mean a family $\mathcal{F}$ of pairwise disjoint nonempty subsets of a set $X$ such that $X = \bigsqcup \mathcal{F}$. 
Note that foliations and the set of orbits of flows are decompositions of connected elements. 
Since connectivity is not required, the set of orbits of group-actions are also decompositions. 
For any $x \in X$, denote by $\F(x)$ the element of a decomposition $\mathcal{F}$ containing $x$. 
Denote by $X/\mathcal{F}$ the quotient space, called the decomposition space. 
In other words, the decomposition space $X/\mathcal{F}$ is a quotient space $X/\sim_{\F}$ defined by $x \sim_{\F} y$ if $\mathcal{F}(x) = \mathcal{F}(y)$. 

\subsubsection{Saturations and their class spaces}
Let $\mathcal{F}$ be a decomposition of a topological space $X$. 
A subset $A \subseteq X$ is invariant (or saturated) if $A = \F (A)$. 
A nonempty invariant subset $A$ is minimal if $\overline{\F(x)} =A$ for any $x \in A$. 
For a subset $A$ of $X$, write the saturation $\F(A)  = \mathrm{Sat}_{\F}(A) := \bigcup_{x \in A} {\F}(x)$. 
The (element) class $\hat{\F}(x)$ of an element $\F(x)$ is defined by $\hat{\F}(x) := \{ y \in X \mid  \overline{\F(y)} = \overline{\F(x)} \} \subseteq \overline{\F(x)}$. 
Then the set $\hat{\mathcal{F}} := \{ \hat{\F}(x) \mid x \in X\}$ of classes is a decomposition, called the class decomposition.  
Denote by $X/\hat{\mathcal{F}}$ the quotient space, called the (element) class space. 
In other words, the class space $X/\hat{\mathcal{F}}$ is a quotient space $X/\sim_{\hat{\F}}$ defined by $x \sim_{\hat{\F}} y$ if $\hat{\mathcal{F}}(x) = \hat{\mathcal{F}}(y)$ (i.e. $\overline{\F(x)} = \overline{\F(y)}$). 
Note that the class $\hat{\F}(x)$ of an element $\F(x)$ is minimal if and only if $\hat{\F}(x) =  \overline{\F(x)}$. 
Moreover, the class space $X/{\hat{\F}}$ is the $T_0$-tification (or Kolmogorov quotient) of the decomposition space $X/\mathcal{F}$. 

\subsubsection{Pointwise almost periodic decompositions}
A decomposition is pointwise almost periodic if the set of all closures of elements of $\mathcal{F}$ also is a decomposition. 
Then we observe the following characterization of pointwise almost periodicity. 
\begin{lemma}
The following statements are equivalent for a decomposition $\mathcal{F}$ on a topological space $X$: 

$(1)$ The decomposition $\mathcal{F}$ is pointwise almost periodic. 

$(2)$ The class space $X/\hat{\mathcal{F}}$ is $T_1$.  

$(3)$ The closure of any element of $\mathcal{F}$ is minimal. 
\\
In any case, the closure of an element corresponds to the class $(\mathrm{i.e.}$ $\hat{\mathcal{F}}(x) = \overline{\F(x)}$ for any point $x \in X)$.  
\end{lemma}
\begin{proof}
Suppose that $\mathcal{F}$ is pointwise almost periodic. 
Then any pair of closures of elements either correspond to each other or is disjoint. 
This means that $\hat{\mathcal{F}}(x) = \overline{\F(x)}$ for any point $x \in X$. 
Therefore $X/\hat{\mathcal{F}}$ is $T_1$. 
Suppose that $X/\hat{\mathcal{F}}$ is $T_1$. 
Then 
$\overline{\{ \hat{L} \}} = \{ \hat{L} \}$ in $X/\hat{\mathcal{F}}$. 
By definition of quotient topology, we have $\overline{L} \subseteq \overline{\hat{L}} = \hat{L}$ in $X$ and so $\overline{L} = \hat{L}$. 
Therefore the closure of any element of $\mathcal{F}$ is minimal. 
Suppose that the closure of any element of $\mathcal{F}$ is minimal. 
This means that $\hat{\mathcal{F}}(x) = \overline{\F(x)}$ for any point $x \in X$, and so that $\mathcal{F}$ is  pointwise almost periodic. 
\end{proof}

We recall nets as follows to define ``characteristic 0''. 
\subsubsection{Orders and nets}
A binary relation $\leq$ on a set $P$ is a pre-order if it is reflexive (i.e. $a \leq a$ for any $a \in P$) and transitive (i.e. $a \leq c$ for any $a, b, c \in P$ with $a \leq b$ and $b \leq c$).
A pre-order $\leq$ on $X$ is a partial order if it is antisymmetric (i.e. $a = b$ for any $a,b \in P$ with $a \leq b$ and $b \leq a$).
A partial order $\leq$ is a total order if each elements are comparable (i.e. either $x \leq y$ or $y \leq x$ for any points $x,y$). 
A pre-order set is a set with a pre-order. 
A pre-order set is a directed set if each finite subset of it has an upper bound. 
A mapping $f \colon D \to X$ is a net if $D$ is a directed set and $X$ is a topological space. 

\subsubsection{Concepts of decompositions}
Recall a point $x$ in a topological $X$ is characteristic 0 \cite{K} if $\Fh (x) =D(x)$ for any $x \in X$, where $D(x)$ is its (bilateral) prolongation defined as follows: $D(x) = \{ y \in X \mid y_{\alpha} \in  \F(x_{\alpha}), y_{\alpha} \to y,  \text{ and } x_{\alpha} \to x \text{ for some nets } (y_{\alpha}),  (x_{\alpha}) \subseteq X \}$. 
A decomposition is characteristic 0 if so is each point of it. 
A decomposition $\mathcal{F}$ on a topological space $X$ is upper semi-continuous (or usc) if each element of $\mathcal{F}$ is both closed and compact and for any $L \in \mathcal{F}$ and for any open neighborhood $U$ of $L$ there is an invariant open neighborhood $V$ of $L$ contained in $U$. 
We will show that upper semi-continuity and $L$-stability coincide if $X$ is compact Hausdorff. 
Here, the decomposition $\mathcal{F}$ is $L$-stable if for any neighborhood $U$ of each element class $\hat{L} \in \hat{\mathcal{F}}$, there is an open ${\F}$-invariant neighborhood of $\hat{L}$ contained in $U$. 
The $L$-stability is an analogous concept of Lyapunov stability. 
A decomposition $\mathcal{F}$ on a topological space is weakly almost periodic in the sense of Gottschalk W. H. \cite{EN} if it is pointwise almost periodic and the saturation $\hat{\mathcal{F}}(A)$ of closures of elements for any closed subset $A$ is closed. 
Since any closures of elements of a pointwise almost periodic decomposition corresponds to the element classes, notice that a pointwise almost periodic decomposition $\mathcal{F}$ on a topological space $X$ is weakly almost periodic if and only if the quotient map $p \colon  X \to X/\hat{\mathcal{F}}$ is closed (i.e. the saturation $\hat{\mathcal{F}}(A)$ for any closed subset $A$ is closed). 
A decomposition $\mathcal{F}$ on a topological space is $R$-closed if the element closure relation $R := \{ (x, y) \mid y \in \overline{\mathcal{F}(x)} \}$ is closed. 
By definition of $R$, notice that $\bigsqcup_{y \in X} (\F(y) \times \{ y \}) \subset R$. 
It is known that an $R$-closed decomposition of a topological space $X$ is pointwise almost periodic (cf. \cite[Corollary~1.4]{Y}). 
Moreover, we have the following observation. 

\begin{lemma}
The element closure relation $R$ of a decomposition on a topological space is symmetric if and only if the decomposition is pointwise almost periodic. 
\end{lemma}

\begin{proof}
Let $\mathcal{F}$ be a decomposition on a topological space $X$. 
Suppose that $R$ is symmetric. 
Fix any $x, y \in X$. 
Then $x \in \overline{\F(y)}$ if and only if $y \in \overline{\F(x)}$.
This implies that if $x \in \overline{\F(y)}$ then $\F(y) \subseteq \overline{\F(x)}$ and so $\overline{\F(y)} \subseteq \overline{\F(x)}$. 
Fix a point $x \in \overline{\F(y)}$. 
Then $\overline{\F(y)} \subseteq \overline{\F(x)}$. 
Since $y \in \overline{\F(x)}$, we have $\overline{\F(x)} \subseteq \overline{\F(y)}$. 
This means that $\overline{\F(x)} = \overline{\F(y)}$. 
Thus $\overline{\mathcal{F}(y)} = \hat{\F}(y)$ and so that the decomposition $\mathcal{F}$ is pointwise almost periodic. 
Conversely, suppose that $\mathcal{F}$ is pointwise almost periodic. 
Then each element is minimal and so $\hat{\mathcal{F}}(x) = \overline{\mathcal{F}(x)}$ for any $x \in X$. 
Therefore, for any $x, y \in X$, we have that $y \in \overline{\mathcal{F}(x)}$ if and only if $x \in \overline{\mathcal{F}(y)}$. 
This means that $R$ is symmetric. 
\end{proof}

%
%
%

\subsection{Foliations}

A foliation is pointwise almost periodic (resp. $R$-closed) if it is pointwise almost periodic (resp. $R$-closed) as a decomposition. 
By the definition of almost periodicity, compact foliations on compact manifolds are pointwise almost periodic.  
In \cite{Y}, it is shown that the set of compact Hausdorff foliations and of minimal foliations on compact manifolds is a proper subset of the set of $R$-closed foliations. 

For a non-negative integer $k$, a foliation $\F$ is codimension-$k$-like if all but finitely many leaf classes of $\F$ are codimension $k$ submanifolds without boundaries and the finite exceptions are  leaf classes which are closed subsets and whose codimensions are more than $k$.   
%
%
The following remark says that the existence of codimension-$k$-like $R$-closed foliations is wider than one of codimension $k$ $R$-closed foliations.

\begin{rem}
There are no codimension one $R$-closed foliations on $\mathbb{S}^3$.  
On the other hand, there are codimension-one-like $R$-closed foliations on $\mathbb{S}^3$. 
\end{rem}

\begin{proof}
By the $C^0$ Novikov Compact Leaf Theorem~\cite{So}, any codimension one foliation on $\mathbb{S}^3$ has Reeb components (i.e. solid tori each of which is foliated by a torus and planes). 
Since any codimension one $R$-closed foliations on a three closed $3$-manifold have no Reeb components, there are no codimension one $R$-closed foliations on $\mathbb{S}^3$. 

On the other hand, consider a minimal irrational flow $v_1$ on a torus $\T^2$, which induces a Kronecker foliation. 
Let $v_2$ be a flow on $\T^2 \times [0,1]$ defined by $v_2(p, s) := (v_1(p), s)$. 
Define an equivalent relation $\sim$ as follows:  $(x,y,s) \sim (x',y',s')$ if either both $s = s' = 0$ and $x = x'$ or both $s = s' = 1$ and $y = y'$. 
Then the quotient space $(\T^2 \times [0,1])/\sim$ is a sphere $\mathbb{S}^3$ such that the induced flow $v$ has exactly two periodic orbits $[\T^2 \times \{ 0 \}]$ and $[\T^2 \times \{ 1 \}]$ and each orbit class except two periodic orbits is a torus.  
This means that the induced foliation $\F_v$ is a codimension-one-like $R$-closed foliations on $\mathbb{S}^3$. 
\end{proof}

\subsection{Group-actions}
An action of a group $G$ is a homomorphism $\varphi_G: G \to \mathrm{Hemeo}(X)$, where $\mathrm{Hemeo}(X)$ is the group of homeomorphism on $X$. 
By a group-action, we mean a continuous action of a topological group $G$ on a topological space $X$. 
We identify, by abuse of notation, groups with their group-actions. 
For a group-action $G$, denote by $\mathcal{F} _G$ the set of orbits of $G$. 
Then $\mathcal{F} _G$ is a decomposition. 
The orbit class $\hat{G}(x)$ of a point $x \in X$ is the subset consisting of points whose orbit closures correspond to the orbit closure of $x$ (i.e. $\hat{G}(x) := \{ y \in X \mid \overline{G(x)} = \overline{G(y)} \}$).  
For a group-action $G$, the class space is the set of orbit classes and is denoted by $M/\hat{G}$ (i.e. $M/\hat{G} := M/\hat{\mathcal{F}}_G$). 
A group-action is pointwise almost periodic if any orbit closure is minimal. 
Note that a group-action $G$ is pointwise almost periodic if and only if so is the decomposition $\mathcal{F} _G$. 
A group-action $G$ is $R$-closed if the decomposition $\mathcal{F} _G$ is $R$-closed. 
Therefore $G$ is $R$-closed if and only if $R := \{ (x, y) \mid y \in \overline{G(x)} \}$ is closed. 
Notice that the relation $R$ for a group-action is called the orbit closure relation. 

For a non-negative integer $k$, a group-action $G$ is codimension-$k$-like if all but finitely many orbit closures of $\F_G$ are codimension $k$ connected submanifolds without boundaries and the finite exceptions are connected closed subsets each of whose codimension is more than $k$.  
For a property P, a group is virtually P if there is a finite index subgroup has the property P. 
Similarly, we call that a group-action $G$ is a virtually codimension-$k$-like group-action by a finite index normal subgroup if there is a finite index normal subgroup whose action is codimension-$k$-like. 


\subsubsection{Flows}
By a flow, we mean a continuous $\R$-action on a topological space $X$. 
Therefore a fow has a property P if so is the flow as a group-action. 
We say that a flow is trivial if it is identical or minimal. 
Notice that the minimal sets of an identical flow is a singleton and the orbit class space is the whole space $M$, and that the minimal set of a minimal flow is the whole space $M$ and the orbit class space is a singleton. 
For a flow $v$  on a $3$-manifold, a minimal set $\mathcal{M}$ of $v$ is called a suspension of a circloid for $v$ if there are an open annulus $\mathbb{A}$ transverse to  $v$ and a circloid $C$ in $\mathbb{A}$ whose saturation $\F_v(C)$ is $\mathcal{M}$ with $C = \mathcal{M} \cap \mathbb{A}$, where $\F_v(C) = \bigcup_{x \in C} O_v(x)$.  

For a flow $v$ on a manifold $M$, denote by $\F_v$ the decomposition of $M$ into orbits of $v$. 
The flow $v$ is $R$-closed if and only if the decomposition $\F_v$ is $R$-closed.
Note $\hat{\mathcal{F}}_{v} = \{ \hat{O}_v(x) \mid x \in M \}$, where the orbit class $\hat{O}_v(x)$ is the subset consisting of elements whose orbit closures correspond to the orbit closure of $x$ (i.e. $\hat{O}_v(x) = \{ y \in M \mid \overline{O_v(x)} = \overline{O_v(y)} \}$).  
For a flow $v$, write $M/v :=M/\mathcal{F}_{v}$ (resp. $M/\hat{v} :=M/\hat{\mathcal{F}}_{v}$), called the orbit space (resp. (orbit) class space) of $v$. 
The following remark implies the non-completeness of orbit spaces. 

\begin{rem} 
There are pairs of flows on a torus which are not topologically equivalent but whose orbit $($class{\rm)} spaces are homeomorphic. 
\end{rem}

\begin{proof}
Consider two minimal linear flows on the flat torus $(\R/\Z)^2$ generated by vector fields $v, v'$ such that $g(v, \cdot ) = [dx_1 + \sqrt{2} dx_2],  g(v', \cdot ) = [dx_1 + \pi dx_2] \in H^1((\R/\Z)^2, \R) \cong \mathbb{R}^2$ respectively, where $g$ is the induced flat Riemann metric from the standard Euclidean metric on $\R^2$.  
Then the orbit spaces are indiscrete spaces with cardinalities of the continuum and the orbit class spaces are singletons. 
On the other hand, assume that $v'$ is topologically equivalent to $v$ via a homeomorphism $h \colon  (\R/\Z)^2 \to (\R/\Z)^2$. 
Then $h^*(g(v, \cdot )) = g(v', \cdot ) \in H^1((\R/\Z)^2, \R)$, where $h^*$ is the pullback along $h$ on $H^1((\R/\Z)^2, \R)$. 
Since the group of the induced mappings on $H^1((\R/\Z)^2, \R)$ by orientation-preserving toral homeomorphisms can be identified with $\operatorname{SL}(2,\mathbb{Z})$, there is a matrix $A \in \operatorname{SL}(2,\mathbb{Z})$ such that $A(1, \sqrt{2})^T = (1, \pi)^T$ if $h$ is orientation-preserving and that $A(1, \sqrt{2})^T = (-1, \pi)^T$ if $h$ is orientation-reserving. 
This means that there are two integers $n,m \in \Z$ such that $\pi = n + m \sqrt{2}$, which contradicts that $\pi$ is transcendental and $\sqrt{2}$ is algebraic. 
\end{proof}

\subsubsection{Homemorphisms}
Note that a homeomorphism can be seen as a continuous $\Z$-action and so a group-action. 
Therefore a homeomorphism $f$ has a property P if so is the group-action generated by $f$. 
For a homeomorphism $f$ on a topological space $X$, denote by $\F_f$ the decomposition of $X$ into orbits of (the group-action generated by) $f$. 

\subsection{Examples of (non-)codimension-$k$-like (non-)$R$-closed flows}

We state a non-$R$-closed example. 

\begin{example}[Non-$R$-closed cases]
A Morse-Smale flow on a closed $n$-manifold is codimension-$(n-1)$-like but not $R$-closed.
\end{example}

\subsubsection{Examples of codimension-one-like flows}
We state codimension-one-like (non-)$R$-closed examples as follows. 

\begin{example}[Codimension-one-like non-$R$-closed cases]
A Denjoy foliation on a torus is codimension-one$($-like{\rm)} but not $R$-closed. 
Moreover, a Morse-Smale flow on a closed surface is codimension-one-like but not $R$-closed.
\end{example}

\begin{example}[Codimension-one-like $R$-closed flows on surfaces]
A spherical flow which is a non-trivial rotation with respect to an axis is not codimension-one but codimension-one-like and $R$-closed, and the orbit class space is a closed interval. 
Moreover, there are codimension-one-like $R$-closed spherical homeomorphisms $(\mathrm{resp.}$ diffeomorphisms{\rm)} with minimal sets which are not circles but circloids \cite{R,W} $(\mathrm{resp.}$ \cite{He,N}$)$. 

Note that a toral flow $v$ defined by $v_t(x,y) = (x + t \sin (2 \pi y), y) \in (\R/\Z)^2$ is  pointwise almost periodic but neither non-codimension-$k$-like for any $k =1,2$ nor $R$-closed. 
\end{example}

\begin{example}[Codimension-one-like $R$-closed flows on $3$-manifolds]
For an irrational number $\alpha \in \R - \Q$, let $f$ be a spherical diffeomorphism defined by 
$$f(\sqrt{1 - z^2} \cos \theta, \sqrt{1 - z^2} \sin \theta ,z) := (\sqrt{1 - z^2} \cos (\theta + \alpha), \sqrt{1 - z^2} \sin (\theta + \alpha) ,z) \in \mathbb{S}^2$$ $(\mathrm{i.e.}$ an irrational rotation with respect to the $z$-axis on the sphere $\mathbb{S}^2 )$, where $\mathbb{S}^2 := \{ (x,y,z) \in \R^3 \mid x^2 + y^2 + z^2 = 1 \}$.  
Denote by $v_f$ the suspension flow on the mapping torus $(\mathbb{S}^2 \times \R)/(x, r) \sim (f(x), r+1)$, which is a closed $3$-manifold. 
Then the $v_f$ is a codimension-one-like $R$-closed flow on the mapping torus $\mathbb{S}^2 \times \mathbb{S}^1$ and so the orbit class space is a closed interval. 
\end{example}

\subsubsection{Examples of codimension-two-like flows}
We also state codimension-two-like $R$-closed examples as follows. 

\begin{example}[Suspension flow of a surface diffeomorphism]
The suspension flow of a periodic spherical diffeomorphism $(\mathrm{e.g.}$ a rational rotation{\rm)} is a codimension-two-like $R$-closed flow on the mapping torus  $\mathbb{S}^2 \times \mathbb{S}^1$,  and the Hopf fibration on the $3$-sphere $\mathbb{S}^3$ is a codimension-two-like $R$-closed foliation whose class space is a sphere. 
Moreover each Seifert fibration is a codimension-two $R$-closed foliation. 

Note that {\rm(}the suspension flow of {\rm)} a spherical diffeomorphism $f \colon  \mathbb{S}^2 \to \mathbb{S}^2$ defined by 
$f(\sqrt{1 - z^2} \cos \theta, \sqrt{1 - z^2} \sin \theta ,z) := (\sqrt{1 - z^2} \cos (\theta + z), \sqrt{1 - z^2} \sin (\theta + z) ,z)$ is pointwise almost periodic but neither codimension-$k$-like nor $R$-closed. 
\end{example}

Recall that the topological suspension $SX$ of a topological space $X$ is the quotient space defined by 
${\displaystyle SX := (X\times I)/\sim}$, where ${\displaystyle  (x_{1},s_1)\sim (x_{2},s_2)}$ if $s_1 = s_2 \in \{ 0, 1\}$. 

\begin{example}[Topological suspension of a surface homeomorphism]
Since the $3$-sphere $\mathbb{S}^3$ is a topological suspension of the sphere $\mathbb{S}^2$, the induced homeomorphism $Sf$ $(\mathrm{i.e.}$ $Sf([x,s]) := [f(x), s])$ on the topological suspension $\mathbb{S}^3 = S\mathbb{S}^2$ of a codimension-one-like $R$-closed spherical homeomorphism $f$ is codimension-two-like and $R$-closed, and the class space is either a sphere or a closed disk with two corners. 
\end{example}

\section{Characterizations of $R$-closedness and upper semi-continuity}
In this section, we characterize $R$-closedness of decompositions and upper semi-continuity for class decompositions. 
First, we show the following relations. 

\begin{lemma}\label{lem:001}
Let $\mathcal{F}$ be a decomposition of a Hausdorff space $X$.  
If the class decomposition $\hat{\mathcal{F}}$ is usc, then $\mathcal{F}$ is $R$-closed and the class space $X/\hat{\mathcal{F}}$ is Hausdorff. 
Conversely, if the closures of elements of $\F$ are compact and $\mathcal{F}$ is $R$-closed, 
then $\mathcal{F}$ is weakly almost periodic and $\hat{\mathcal{F}}$ is usc. 
\end{lemma}

\begin{proof}
Suppose that $\hat{\mathcal{F}}$ is usc.  
By  \cite[Proposition~1.2.1, p.13]{D}, we have that $X/\hat{\mathcal{F}}$ is Hausdorff. 
By  \cite[Lemma~2.3]{Y}, we obtain that $\mathcal{F}$ is $R$-closed. 
Conversely, suppose that the closures of elements of $\F$ are compact and that $\mathcal{F}$ is $R$-closed. 
We claim that \cite[Lemma~1.6]{Y} is true when the Hausdorff separation axiom for $X$ is replaced by the compact orbit closure condition. 
Indeed, the proof of \cite[Lemma~1.5]{Y} is true when the Hausdorff separation axiom for $X$ is replaced by the compact orbit closure condition. 
\cite[Lemma~1.6]{Y} is true when the Hausdorff separation axiom for $X$ is replaced by the compact orbit closure condition. 

Therefore $\mathcal{F}$ is weakly almost periodic and so that the quotient map $p \colon  X \to X/\hat{\mathcal{F}}$ is closed. 
Since each element of $\hat{\mathcal{F}}$ is compact, from \cite[Proposition~1.1.1, p.8]{D}, we have that $\hat{\mathcal{F}}$ is usc.
\end{proof}

Notice that the compact condition is necessary. 
Indeed, the $R$-closed decomposition $\F = \{ X \}$ for a non-compact Hausdorff space $X$ is weakly almost periodic but $\hat{\mathcal{F}} = \F$ is not usc. 
We show a correspondence for invariances. 

\begin{lemma}\label{lem:bdry_sat}
Let $\mathcal{F}$ be a decomposition of a topological space $X$. 
A closed subset of $X$ is $\mathcal{F}$-invariant if and only if it is $\hat{\F}$-invariant. 
Moreover, we obtain $F = \F (F) = \hat{\F}(F) = \bigcup_{x \in F} \overline{\F(x)}$ for any closed invariant subset $F \subseteq X$. 
\end{lemma}

\begin{proof}
Let $F$ be a closed subset of $X$. 
Suppose that $F$ is $\mathcal{F}$-invariant. 
Since $F$ is closed, we have $\F (F) = F = \bigcup_{x \in F} \overline{\F(x)} = \bigcup_{x \in F} \hat{\F}(x) = \hat{\F}(F)$. 
Since $\hat{\F}$-invariant subsets are $\F$-invariant, the converse holds. 
\end{proof}

We show that each open $\mathcal{F}$-invariant subset for a decomposition $\F$ is $\hat{\mathcal{F}}$-invariant.  

\begin{lemma}\label{lem02}
Let $\mathcal{F}$ be a decomposition of a topological space.  
Then an open subset is $\mathcal{F}$-invariant if and only if it is $\hat{\mathcal{F}}$-invariant.  
Moreover, we obtain $U = \F(U) = \hat{\F}(U)$ for any invariant open subset $U \subseteq X$. 
\end{lemma}

\begin{proof} 
Let $U$ be an open subset of a topological space $X$.  
Suppose that $U$ is ${\F}$-invariant. 
Then the complement $X - U$ is a closed ${\F}$-invariant subset of $X$. 
Lemma~\ref{lem:bdry_sat} implies that the closed subset $X - U$ is $\hat{\F}$-invariant and so $U$ is open $\hat{\F}$-invariant. 
Conversely, suppose that $U$ is $\hat{\F}$-invariant. 
Since any $\hat{\F}$-invariant subsets are ${\F}$-invariant, the subset $U$ is ${\F}$-invariant. 
\end{proof}

Notice that ${\F}$-invariant subsets are not $\hat{\F}$-invariant in general. 
In fact, any orbits for an irrational rotation $f \colon \mathbb{S}^1 \to \mathbb{S}^1$ on a circle are ${\F}_f$-invariant but not $\hat{\F}_f$-invariant. 
Using the previous lemma, we have the following statement. 

\begin{lemma}\label{lem:001a}
Let $\mathcal{F}$ be a pointwise almost periodic decomposition of a normal space $X$.  
Suppose that,  for any open neighborhood $W$ of each element $L \in \mathcal{F}$, there is an invariant open neighborhood of $L$ contained in $W$. 
Then $X/\hat{\mathcal{F}}$ is Hausdorff.  
\end{lemma}

\begin{proof}
Suppose that, for any open neighborhood $W$ of each element $L \in \mathcal{F}$, there is an open $\mathcal{F}$-invariant neighborhood of $L$ contained in $W$. 

We claim that,  for any open neighborhood $U$ of each element $\hat{L} \in \hat{\mathcal{F}}$, there is an open $\hat{\mathcal{F}}$-invariant neighborhood of $\hat{L}$ contained in $U$. 
Indeed, fix any open neighborhood $U$ of each element $\hat{L} \in \hat{\mathcal{F}}$. 
The hypothesis implies that, for any $L' \in \mathcal{F}$ contained in $\hat{L}$, since $U$ is also an open neighborhood of $L'$,  there is an open $\mathcal{F}$-invariant neighborhood $U_{L'}$ of $L'$ contained in $U$.
By Lemma~\ref{lem02}, the union  $\bigcup \{ U_{L'} \mid L' \in \mathcal{F}, L' \subseteq \hat{L} \}$ is desired. 

Fix any elements $\hat{L} \neq \hat{L}' \in \hat{\F}$ if exists. 
Since $\F$ is pointwise almost periodic, we obtain $\hat{L} = \overline{L}$ and $\hat{L}' = \overline{L'}$. 
Since $X$ is normal, there are disjoint open neighborhoods $U_0$ and $V_0$ of $\hat{L}$ and $\hat{L}'$ respectively. 
By the claim,  there are open $\hat{\F}$-invariant neighborhoods $U_1 \subseteq U_0$ and $V_1 \subseteq V_0$ of $\hat{L}$ and $\hat{L'}$ respectively such that $U_1 \cap V_1 = \emptyset$.  
\end{proof}

Notice that the converse is not true. 
Indeed, consider $X = \mathbb{T}^2 = \R/\Z \times \R/\Z$, $\alpha \in \R - \Q$,  and homeomorphism $f$ of $X$ by $f(x, y) := (x + \alpha, y)$. 
Then the class space $\mathbb{T}^2/\hat{\F_v}$ is a circle and so is Hausdorff. 
Then the saturation $\F_v(A)$ of  $A := \{ (x, y) \mid y \in [0, 1/2],  x \in [1/2, 1/2 + y] \}$ contains $ \R/\Z \times (0, 1/2]$ and so the open neighborhood $\mathbb{T}^2 - A$ of the orbit of $0$ contains no invariant open neighborhood. 
We have the following statement.

\begin{corollary}\label{cor:equiv_R}
An $L$-stable decomposition $\mathcal{F}$ of a $T_1$ space $X$ satisfies that $X/\hat{\mathcal{F}}$ is Hausdorff. 
Moreover, if $X$ is $T_3$, then $\mathcal{F}$ is $R$-closed. 
\end{corollary}

\begin{proof}
Let $\mathcal{F}$ be an $L$-stable decomposition of a $T_1$ space $X$.
\cite[Lemma~3.2]{Y} implies that $\mathcal{F}$ is pointwise almost periodic. 
Applying Lemma~\ref{lem:001a} to the class decomposition $\hat{\mathcal{F}}$, 
the class space $X/\hat{\mathcal{F}}$ is Hausdorff.  
Suppose that $X$ is $T_3$. 
By \cite[Lemma~2.2]{Y}, the decomposition $\mathcal{F}$ is $R$-closed. 
\end{proof}

The weakly almost periodicity implies the $L$-stability. 

\begin{lemma}\label{lem:wap_L_stable}
A weakly almost periodic decomposition of a topological space is $L$-stable. 
Moreover, if the closure of every element is compact, then the class decomposition is usc. 
\end{lemma}

\begin{proof}
Let $\mathcal{F}$ be a weakly almost periodic decomposition of a topological space $X$ and $L \in \F$ an element. 
Fix any open neighbourhood $U$ of the class $\hat{L} \in \hat{\F}$. 
Since the complement $X - U$ is closed, the weakly almost periodicity implies that the saturation $\hat{\F}(X -U)$ is closed. 
By $\hat{L} \cap (X - U) = \emptyset$, we have $\hat{L} \cap \hat{\F}(X -U) = \emptyset$. 
Therefore the complement $X - \hat{\F}(X -U) \subseteq U$ is a $\hat{\F}$-invariant open neighborhood of $\hat{L}$ and so is $\F$-invariant.
\end{proof}

\subsection{Proof of Theorem~\ref{th:usc}}
Let $\mathcal{F}$ be a decomposition of a $T_4$ space $X$.  
Since each of the conditions (1) and (2) implies that $\F$ is pointwise almost periodic, by  \cite[Lemma~2.2]{Y}, the conditions (1) and (2) are equivalent. 
 \cite[Corollary~3.1]{Y} implies that the conditions (1) and (3) are equivalent. 
By Corollary~\ref{cor:equiv_R}, the condition $(5)$ implies the $R$-closedness of $\F$. 
From Lemma~\ref{lem:001}, the condition $(4)$ leads to the $R$-closedness of $\F$. 
By Lemma~\ref{lem:wap_L_stable}, the condition $(6)$ implies $L$-stability of $\F$ and so $R$-closedness of $\F$. 

Suppose that $X$ is compact. 
 \cite[Corollary~3.4]{Y} implies that the conditions (1) and (5) are equivalent. 
By Lemma~\ref{lem:001}, the conditions (1) and (4) are equivalent and the condition (1) implies the condition (6). 
This completes a proof of Theorem~\ref{th:usc}.

\subsection{Closures and saturations of $R$-closed decompositions}

We have the following equivalence. 

\begin{lemma}\label{lem:cl_sat}
Let $\mathcal{F}$ be a decomposition of a compact Hausdorff space $X$.  
Suppose that the closure of any invariant subset of $X$ is invariant. 
The following conditions are equivalent: 
\\
$(1)$ The decomposition $\mathcal{F}$ is $R$-closed.  
\\
$(2)$ For any subset $A \subseteq X$, we have $\overline{\mathcal{F}(A)} =  \hat{\mathcal{F}}(\overline{A})$.  

In any case, we have $\overline{\mathcal{F}(A)} =  \hat{\mathcal{F}}(\overline{A}) =  \mathcal{F}(\overline{A})$ for any subset $A \subseteq X$. 
\end{lemma}

\begin{proof} 
Suppose that the closure of any invariant subset of $X$ is invariant. 
For any subset $A \subseteq X$, since $\overline{A} \subseteq \overline{\mathcal{F}(A)}$, 
we obtain $\mathcal{F}(A) \subseteq \hat{\mathcal{F}}(\overline{A}) \subseteq \overline{\mathcal{F}(A)}$. 
Suppose that $\mathcal{F}$ is $R$-closed.  
Fix a subset $A \subseteq X$. 
From Lemma~\ref{lem02}, Theorem~\ref{th:usc} implies that the saturation $\mathcal{F}(\overline{A}) = \hat{\mathcal{F}}(\overline{A})$ of the closed subset $\overline{A}$ is closed and so $\overline{\mathcal{F}(A)} =  \hat{\mathcal{F}}(\overline{A})$.  
Conversely, suppose that $\overline{\mathcal{F}(A)} =  \hat{\mathcal{F}}(\overline{A})$ for any subset $A \subseteq X$. 
Fix a subset $A \subseteq X$. 
Lemma~\ref{lem02} implies that $\overline{\mathcal{F}(A)} =  \hat{\mathcal{F}}(\overline{A}) = \mathcal{F}(\overline{A})$ and so that the $\hat{\F}$-saturation of any closed subset of $X$ is closed. 
By Theorem~\ref{th:usc}, the decomposition $\mathcal{F}$ is $R$-closed.  
\end{proof} 

The $R$-closedness for group-actions implies that the closure of the saturation corresponds to the saturation of the closure. 

\begin{proposition}\label{corollary:cl_sat}
Let $G$ be a group-action on a compact Hausdorff space $X$.  
The following conditions are equivalent: 
\\
$(1)$ The group-action $G$ is $R$-closed.  
\\
$(2)$ For any subset $A \subseteq X$, we have $\overline{G(A)} = \hat{G}(\overline{A})$.  

In any case, we have $\overline{G(A)} = \hat{G}(\overline{A}) = G(\overline{A})$ for any subset $A \subseteq X$. 
\end{proposition}

\section{Inherited properties of $R$-closed group-actions}
In this section, we show the inherited properties of $R$-closed group-actions. 
To show such properties, we show the following statements. 

\begin{lemma}\label{lem:closed}
Let $\pi \colon  G \times X \to X$ be a group-action on a topological space $X$. 
For a compact subset $K \subseteq G$ and for a closed subset $C \subseteq X$, we have that $K \cdot C$ is closed. 
\end{lemma}

\begin{proof} 
Fix a point $x \in X - K \cdot C$. 
Since $x \neq k \cdot c$ for any $c \in C$ and $k \in K$, we have $k^{-1} \cdot x \neq c$ for any $c \in C$ and $k \in K$. 
This means that $(K^{-1} \cdot x) \cap C = \emptyset$. 
Then the complement $X - C$ is an open neighborhood of $K^{-1} \cdot x$ and so $\pi^{-1}(X - C)$ is an open neighborhood of $K^{-1} \times \{ x \}$. 
For any $k \in K^{-1}$, there are neighborhoods $V_k \subseteq G$ and $U_k \subseteq X$ of $k$ and $x$ respectively such that $V_k \times U_k \subseteq \pi^{-1}(X - C)$. 
Since $K^{-1}$ is compact, there is a finite subset $F$ of $K^{-1}$ such that $K^{-1} \subseteq \bigcup_{k \in F} V_k$.  
Then $U := \bigcap_{k \in F} U_k$ is an open neighborhood of $x$ with $K^{-1} \times U \subseteq \bigcup_{k \in F} V_k \times U \subseteq \bigcup_{k \in F} V_k \times U_k \subseteq \pi^{-1}(X - C)$. 
Hence  $K^{-1} \cdot U \subseteq X - C$. 
Therefore $(K^{-1} \cdot U) \cap C = \emptyset$ and so $U \cap (K \cdot C) = \emptyset$. 
This shows that $K \cdot C$ is closed. 
\end{proof}

\begin{lemma}\label{lem:action_closure}
Let $G$ be a group-action on a topological space $X$ and $H$ a subset of $G$. 
Then $\overline{g \cdot (H \cdot x)} = g \cdot \overline{H \cdot x}$ for any $g \in G$. 
\end{lemma}

\begin{proof} 
Fix $g \in G$. 
Since $g$ is a homeomorphism, the element $g$ is a closed mapping. 
By $g \cdot (H \cdot x) \subseteq g \cdot \overline{H \cdot x}$, we have $\overline{g \cdot (H \cdot x)} \subseteq g \cdot \overline{H \cdot x}$. 
Moreover, we obtain $H \cdot x \subseteq g^{-1} \cdot \overline{g \cdot (H \cdot x)} \subseteq  \overline{H \cdot x}$.
Since $g^{-1}$ is a closed mapping, we have $g^{-1} \cdot \overline{g \cdot (H \cdot x)} = \overline{H \cdot x}$ and so $\overline{g \cdot (H \cdot x)} = g \cdot \overline{H \cdot x}$. 
\end{proof}

Recall that a subgroup $H$ of a topological group $G$ is syndetic if there is a compact subset $K$ of $G$ such that $K \cdot H = G$. 

\begin{lemma}\label{lem:synd}
Let $G$ be a group-action on a topological space $X$, $H$ a syndetic subgroup of $G$, and $K \subseteq G$ a compact subset with $K \cdot H = G$. 
If $H$ is $R$-closed, then so is $G$ with $\overline{G \cdot x} = K \cdot \overline{H \cdot x}$ for any $x \in X$. 
\end{lemma}

\begin{proof} 
Let $R_G := \{ (x, y) \mid y \in \overline{G \cdot x} \}$ and $R_H := \{ (x, y) \mid y \in \overline{H \cdot x} \}$. 
Suppose that $H$ is $R$-closed.  
%
Lemma~\ref{lem:action_closure} implies that $\overline{G \cdot x} = {\overline{K \cdot (H \cdot x)}} \supseteq {\bigcup_{g \in K} \overline{g \cdot (H \cdot x)}} = {\bigcup_{g \in K} g \cdot  \overline{(H \cdot x)}} = {K \cdot \overline{H \cdot x}}$. 
By Lemma~\ref{lem:closed}, since $\overline{H \cdot x}$ is closed, we obtain ${K \cdot \overline{H \cdot x}} = \overline{K \cdot \overline{H \cdot x}} \supseteq \overline{K \cdot (H \cdot x)} = \overline{G \cdot x} \supseteq {K \cdot \overline{H \cdot x}}$ and so ${K \cdot \overline{H \cdot x}}= \overline{G \cdot x}$. 
Consider an action of $G$ on $X \times X$ by $g \cdot (x, y) := (x , g^{-1} \cdot y)$.  
Then $K \cdot R_H = \bigcup_{k \in K}  \{ (x, y) \mid k^{-1} y \in \overline{H \cdot x} \}= \bigcup_{k \in K}  \{ (x, y) \mid  y \in k \cdot \overline{H \cdot x} \} =  \{ (x, y) \mid  y \in K \cdot \overline{H \cdot x} \} =  \{ (x, y) \mid  y \in \overline{G \cdot x} \} = R_G$. 
Hence $R_G = K \cdot R_H$. 
Since $R_H$ is closed,  applying Lemma~\ref{lem:closed} to the group-action $G \times (X \times X) \to (X \times X)$,  we have that $R_G = K \cdot R_H$ is closed.   
\end{proof}

\subsection{Proof of Thereom~\ref{main:c}}
%
%
Let $G$ be a group-action on a locally connected compact Hausdorff space $X$ and $H$ a finite index normal subgroup. 
By Lemma~\ref{lem:synd}, the $R$-closedness of $H$ implies one of $G$. 
Conversely, suppose that $G$ is $R$-closed. 
Put $\F := \F_{G}$. 
By  \cite[Corollary~1.4]{Y}, we have that  $G$ is pointwise almost periodic. 
By  \cite[Theorem~2.24]{GH}, we have that $H$ is also pointwise almost periodic.  
By Theorem~\ref{th:usc}, the decomposition $\hat{\F}$ is usc and it suffices to show that $\hat{\F}_H$ is usc. 
Fix any $x \in X$. 
Note that $\F_H (x) = H \cdot x \subseteq G \cdot x = \F (x) $. 
The pointwise almost periodicity implies $\hat{\F}_H (x) = \overline{\F_H (x)} \subseteq \overline{\F(x)} = \hat{\F}(x)$. 
Suppose that $\hat{\F}_H (x) = \hat{\F} (x)$. 
For any open neighborhood $U $ of $\hat{\F} (x) = \hat{\F}_H (x)$, the upper semi-continuity of $\hat{\F}$ implies that there is an $\hat{\F}$-invariant open neighborhood $V$ of $\hat{\F} (x)$ contained in $U$. 
Since $\hat{\F}_H (y) \subseteq \hat{\F} (y)$ for any $y \in X$, we have that $V$ is also an ${\hat{\F}_H}$-invariant open neighborhood  $V$ of $\hat{\F}_H (x)$. 

Thus we may assume that $\hat{\F} (x) \neq \hat{\F}_H (x)$. 
Fix an open neighborhood $U_0$ of $\hat{\F}_H (x)$. 
Let $K$ be a finite subset of $G$ such that $G = K \cdot H = H \cdot K$. 
Put $L := \F_H (x)$ and $K' := \{ k \in K \mid \hat{\F}_H (k \cdot x) \cap \hat{\F}_H (x) = \emptyset \}$. 
The pointwise almost periodicity implies $\hat{L} = \hat{\F}_H (x) = \overline{\F_H(x)} = \overline{H \cdot x}$. 
Denote by $L' :=  \F_H (K' \cdot x)$ be the saturation of $K' \cdot x$. 
Since $K'$ is finite, Lemma~\ref{lem:action_closure} implies $\hat{L}' := \overline{L'} = \overline{H \cdot  (K' \cdot x)} = \overline{K' \cdot  (H \cdot x)} = K' \cdot \overline{H \cdot x} = K' \cdot \hat{L}$ and $\hat{L}' = \overline{\bigcup_{k \in K'} \F_H (k \cdot x)} = \bigcup_{k \in K'} \overline{\F_H (k \cdot x)} = \bigcup_{k \in K'} \hat{\F}_H (k \cdot x)$.  
Then $\hat{L}$ and $\hat{L}'$ are disjoint and closed. 
By Lemma~\ref{lem:synd}, the pointwise almost periodicity of $H$ implies that $\hat{L} \sqcup \hat{L}' = \overline{H \cdot x} \sqcup (K' \cdot \overline{H \cdot x}) = K \cdot \overline{H \cdot x}= \overline{K \cdot (H \cdot x)} = \overline{G \cdot x} = \hat{\F}(x)$. 
Since $X$ is normal, there are disjoint open neighborhoods $U$ and $U'$ of $L$ and $L'$ respectively such that $\overline{U} \subseteq U_0$. 
Then $\hat{\F}(x) = \hat{L} \sqcup \hat{L}' \subseteq U \sqcup U'$. 
Since $\hat{\F}$ is usc, there are open subsets $V \subseteq U$ and $V' \subseteq U'$ such that $V \sqcup V'$ is an $\hat{\F}$-invariant open neighborhood of $\hat{\F}(x) = \hat{L} \sqcup \hat{L}'$. 
Since $\hat{L}$ is compact and $X$ is locally connected, there are finitely many connected open subsets $B_1, \ldots, B_k \subseteq V$ such that $B_i \cap \hat{L} \neq \emptyset$ and $\hat{L} \subseteq B_1 \cup \cdots \cup B_k \subseteq V$. 
Put $B := B_1 \cup \cdots \cup B_k$. 
For any $h \in H$ and $i = 1, \ldots, k$, the $\hat{\F}_H$-invariance of  $\hat{L} = \hat{\F}_H (x)$ implies that $h \cdot B_i \cap \hat{L} \neq \emptyset$ and so $h \cdot B_i \cap V \neq \emptyset$. 
The $\hat{\F}$-invariance of $V \sqcup V'$ implies $h \cdot B_i \subseteq V \sqcup V'$ for any $h \in H$ and $i = 1, \ldots, k$. 
Since $h \cdot B_i$ is connected and $h \cdot B_i \cap V \neq \emptyset$ for any $h \in H$,  we obtain $h \cdot B_i \subseteq V$ for any $h \in H$ and $i = 1, \ldots, k$. 
Therefore $H \cdot B \subseteq V$. 
Since $\hat{\F}_H (B) \subseteq \overline{H \cdot B} \subseteq \overline{V} \subseteq \overline{U} \subseteq U_0$, the neighborhood $\hat{\F}_H (B)$ of $\hat{\F}_H (x)$ is desired. 
This completes a proof of Theorem~\ref{main:c}. 

\subsubsection{A question on Thereom~\ref{main:c}}
The author would like to know whether the local connectivity condition is necessary or not, and whether the finite condition can be replaced by the syndetic condition.

\section{Applications to codimension-one-like or codimension-two-like decompositions} 

In this section, we characterize necessary and sufficient conditions for the class spaces of codimension-$k$-like $(k = 1, 2)$ group-actions and foliations to be $k$-manifolds.

\subsection{Codimensionality of Cantor manifolds}
A separable metrizable space $X$ whose small inductive dimension is $n > 0$ is an $n$-dimensional Cantor manifold \cite{U} if the complement $X - L$ for any closed subset $L$ of $X$ whose small inductive dimension is less than $n - 1$ is connected. 
Hurewicz-Menger-Tumarkin Theorem~in dimension theory says that a connected topological manifold (with or without boundary) is a Cantor manifold (cf. \cite[Corollary~1 and Corollary~2 after Theorem~IV 4, p.48]{HW}). 
Moreover, homogeneous continua are Cantor manifolds (cf. \cite[Theorem~1]{Ku1990}). 
Recall that a topological space $X$ is homogeneous if for each two points $x, y \in X$ there is a homeomorphism $h  \colon  X \to X$ with $h(x)=y$. 
In addition, it is known that any open connected subset of a homogeneous locally compact metrizable space is a Cantor manifold \cite[Theorem~2.1]{Ku}. 
Recall that $X$ is locally compact if any point of $X$ has a compact neighborhood.
Moreover, it is known that any open connected subset of a locally compact locally connected separable metric space $X$ with $H_k (X, X - \{ x \}) = 0$ for any $x \in X$ and any $k < n$ is an $n$-dimensional Cantor manifold, where $H_k(X, X - \{ x \})$ are the singular homology groups of $X$ relative to the subspace $X - \{ x \}$ \cite[Proposition~1.7]{Kr}. 
In particular, any open connected subset of a compact manifold is a Cantor manifold \cite{HM,T,U} (cf.  \cite[Example VI 11, p.93]{HW}).
Recall that $\partial A := \overline{A} - \mathrm{int} A$ for a subset $A$ of a topological space. 
We have the following statement. 

\begin{lemma}\label{lem062} 
Let $M$ be a connected compact manifold, $\F$ an $R$-closed decomposition on $M$, and $U \subseteq M$ nonempty invariant open subset.  
If the codimension of $\F (\partial U)$ is more than one, then $\overline{U} = M$.  
\end{lemma}

\begin{proof} 
The $R$-closedness implies that the saturation $\F (\partial U)$ is closed. 
Hurewicz-Menger-Tumarkin Theorem~(cf. \cite[Theorem~2.1]{Ku}) implies that a connected manifold with or without boundary is a Cantor manifold. 
Since the codimension of $\F (\partial U)$ is more than one, the set differnce $M - \F (\partial U)$ is connected.  
Since $U = \overline{U} \setminus \F (\partial U)$ is closed and open in $M - \F (\partial U)$, the connectivity of $M - \F (\partial U)$ implies that $M  =  U \sqcup \F (\partial U)$. 
Since the codimension of $\F (\partial U)$ is more than one,  for any $x \in \F (\partial U)$ and any \nbd $V$ of $x$, we have $U \cap V = (M - \F (\partial U)) \cap V \neq \emptyset$.  
Therefore $\F (\partial U) \subset \overline{U}$ and so $\overline{U} =  U \sqcup \F (\partial U) = M$. 
\end{proof}

\subsection{Codimension-one-like decompositions}

We state the following statement. 

\begin{lemma}\label{lem:4.1}
Let $\F$ be an $R$-closed decomposition on a compact connected manifold $M$ without codimension zero element classes and $U$ a connected component of the union of codimension one elements of $\hat{\mathcal{F}}$. 
Suppose that each element class is connected. 
If $U$ is nonempty and open such that each element of $\hat{\mathcal{F}}$ in $U$ is a connected submanifold without boundary, then $\overline{U} = M$ and the class space $M/\hat{\mathcal{F}}$ is a closed interval or a circle such that there are at most two elements whose codimensions are more than one.  
\end{lemma}

\begin{proof}
Compactness of $M$ and $R$-closedness of $\F$ imply that each element class is compact. 
By Theorem~\ref{th:usc}, we have that $M/\hat{\mathcal{F}}$ is Hausdorff. 
Let $U$ be the above connected component and $p \colon  M \to M/\hat{\mathcal{F}}$ the canonical projection. 
By  \cite[Theorem~3.3]{D2}, we have that $p(U)$ is a $1$-manifold. 
Suppose $p(U)$ is a circle or a closed interval. 
Since $U$ is open and closed, we have $U = M$.  
Suppose that $U$ is an interval which is not closed. 
Since $M/\hat{\mathcal{F}}$ is Hausdorff, each convergence sequence has one limit and so the saturation $\hat{\F}(\partial U) = \F(\partial U)$ consists of at most two elements whose codimension is more than one.
Lemma~\ref{lem062} implies that $\overline{U} = M$. 
Therefore the closure $\overline{p(U)} =p(\overline{U}) = p(M) = M/\hat{\mathcal{F}}$ is a closed interval. 
\end{proof}

\begin{figure}
\begin{center}
\includegraphics[scale=0.3]{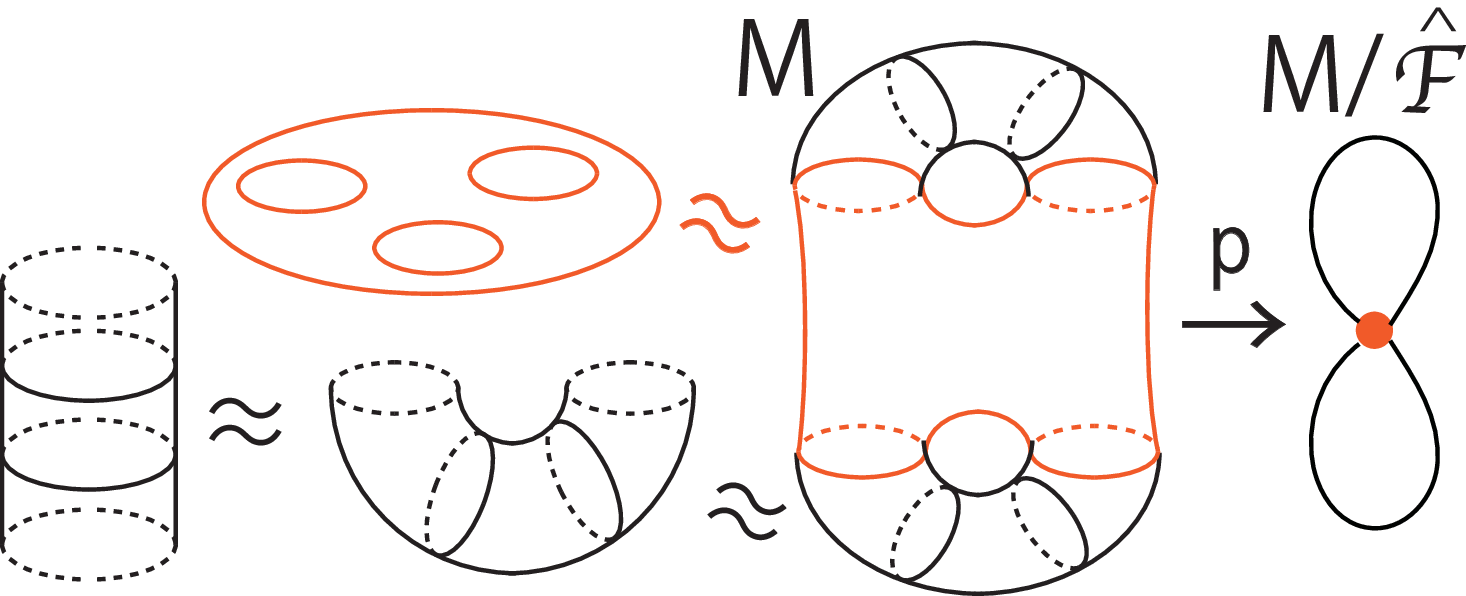}
\end{center} 
\caption{An $R$-closed decomposition on a closed orientable surface with genus $2$ whose class space is a figure eight curve. }\label{fig01}
\end{figure}

The codimension zero condition is necessary, because there are $R$-closed decompositions on compact connected surfaces whose class space is not a manifold. 
Indeed, let $D$ be a compact three-punctured disk with a decomposition $\F_D := \{ D \}$ and $\mathbb{A} = \mathbb{S}^1 \times (0,1)$ an open annulus with a decomposition $\F_{\mathbb{A}} := \{ \mathbb{S}^1 \times \{ y \} \mid y \in (0,1) \}$. 
Gluing $D$ and two copies of $\mathbb{A}$, we obtain a closed surface $S$ of genus $2$ and the induced decomposition $\F$ consists of $D$ and essential circles such that the quotient space $S/\hat{\F}$ is neither an interval nor a circle but a graph which is a figure eight curve (see Figure~\ref{fig01}).  

Recall that each continuous codimension one compact foliation on a compact manifold is $R$-closed. 
Hence the following corollary is a generalization of a fact that the leaf (class) space of a continuous codimension one compact foliation of a compact manifold is a closed interval or a circle.  

\begin{corollary}\label{cor84}
Let $\mathcal{F}$ be an $R$-closed foliation on a compact connected manifold $M$.  
If there is a compact $\hat{\mathcal{F}}$-invariant subset $C$ whose codimension is more than one such that all leaf closures of $\mathcal{F}$ in $M - C$ are codimension one connected submanifolds without boundaries, then $M/\hat{\mathcal{F}}$ is either a closed interval or a circle.  
\end{corollary}

Corollary~\ref{cor84} implies (1) of Theorem~\ref{main:a}.  
Recall that the quotient space of a $1$-manifold by a finite group-action is a $1$-orbifold, which is homeomorphic to a manifold with or without boundary. 
Applying Lemma~\ref{lem:4.1} to the group-action cases, we obtain the following statement. 

\begin{lemma}\label{prop:003}
Let $G$ be an $R$-closed group-action on a compact connected manifold $M$, $H$ a finite index normal subgroup of $G$, and $U$ a connected component of the union of orbit closures of $H$ whose codimension are one.  
If $U$ is nonempty and open such that each orbit closure of $H$ in $U$ is a connected submanifold without boundary, then $M/\hat{\mathcal{F}} _G$ is a closed interval or a circle such that there are at most two elements whose codimensions are more than one.  
\end{lemma}


\begin{proof}
Compactness of $M$ and $R$-closedness of $G$ imply that each orbit class is compact. 
By Theorem~\ref{main:c}, we have that $H$ is also $R$-closed. 
Since $G$ is a non-minimal group-action, so is $H$. 
This implies that there are no codimension zero minimal sets of $H$ because $G/H$ is finite. 
Lemma~\ref{lem:4.1} implies that $\overline{U} = M$ and that the class space $M/\hat{\mathcal{F}}_H$ is either a circle or a closed interval. 
Since $\hat{\mathcal{F}} _H$ consists of orbit closures and since the quotient of a $1$-manifold by a finite group-action is also a $1$-manifolds with or without boundary, the boundary $\partial U$ also consists of at most two orbit closures of $H$ whose codimension are more than one.  
%
Since $G/H$ is a finite group acting $M/\hat{\mathcal{F}}_H$, 
Lemma~\ref{lem:synd} implies that $M/\hat{\mathcal{F}} _G = (M/\hat{\mathcal{F}} _H)/(G/H)$ is compact and so either a closed interval or a circle and so that there are at most two elements whose codimensions are more than one. 
\end{proof}

This lemma implies the following application. 

\begin{proposition}\label{prop:0073}
Let $G$ be an $R$-closed group-action on a compact connected manifold $M$ and $H$ a finite index normal subgroup of $G$. 
Suppose that all but finitely many orbit closures of $H$ are codimension one connected submanifolds without boundaries. 
If each of finite exceptions is closed and its codimension is not one, then $M/\hat{\mathcal{F}} _G$ is a closed interval or a circle.  
\end{proposition}

\begin{proof}
Compactness of $M$ and $R$-closedness of $G$ imply that each orbit class is compact. 
Let $\hat{O}_1, \ldots, \hat{O}_k$ be the finite exceptions and $U := M- \bigcup_{i =1}^k \hat{O}_i$. 
Then $U$ is the union of orbit closures of $H$ which are codimension one connected submanifolds, and is an invariant open submanifold. 
If there is a codimension zero minimal set, then the $R$-closedness implies that it has no boundaries and so equals  $M$, which contradicts to the non-minimality.  
Thus there are no codimension zero minimal sets and so the codimension of each $\hat{O}_i$ is at least two. 
By Lemma~\ref{prop:003}, we are done. 
\end{proof}

Proposition~\ref{prop:0073} implies (1) of Theorem~\ref{main:b}. 

%
%

\subsection{Codimension-two-like decompositions}

Consider the direct system $\{K_a\}$ of compact subsets of a topological space $X$ and inclusion maps such that the interiors of $K_a$ cover $X$.  
There is a corresponding inverse system $\{ \pi_0( X - K_a ) \}$, where $\pi_0(Y)$ denotes the set of connected components of a space $Y$. 
Then the set of ends of $X$, called the end set of $X$, is defined to be the inverse limit of this inverse system.
We have the following sufficient condition for class spaces to be surfaces with corners. 

\begin{lemma}\label{lem:4.2}
Let $\mathcal{F}$ be an $R$-closed decomposition on a compact manifold $M$.  
Suppose that all but finitely many element closures of $\mathcal{F}$ are codimension two connected submanifolds without boundaries. 
If each of finite exceptions is connected and closed and has codimension at least three, then $M/\hat{\mathcal{F}}$ is a surface with corners. 
\end{lemma}

\begin{proof} 
Compactness of $M$ and $R$-closedness of $\F$ imply that each element class is compact. 
By Theorem~\ref{th:usc}, the class decomposition $\hat{\mathcal{F}}$ is usc. 
Let $\hat{L}_1, \ldots, \hat{L}_k$ be all higher codimension elements of $\hat{\mathcal{F}}$. 
Removing higher codimensional elements, let $M' := M - \bigsqcup_{i=1}^k \hat{L}_i$ be the resulting manifold and $\hat{\mathcal{F}} '$ the resulting decomposition of $M'$. 
Then $\hat{\mathcal{F}} '$ consists of codimension two closed connected submanifolds and is usc. 
By \cite[ Theorem~3.12]{D3}, the class space $M'/\hat{\mathcal{F}}'$ is a surface $S'$ with corners.  
Then $(M/\hat{\mathcal{F}} ) - \{\hat{L}_1, \ldots, \hat{L}_k \} \cong M'/\hat{\mathcal{F}} ' = S'$.  

We claim that $S'$ has $k$ ends. 
Indeed, since the finite exceptions $\hat{L}_i$ are connected, there are pairwise disjoint neighborhoods $U_i$ of them.
Since $\hat{\mathcal{F}}$ is usc, there are pairwise disjoint invariant open neighborhoods $V_i \subseteq U_i$ of them. 
Since $W_i - \hat{L}_i$ is connected for any connected neighborhoods $W_i$ of $\hat{L}_i$, each end of $S'$ corresponds to  some $\hat{L}_i$. 
This shows that $S'$ has $k$ ends corresponding to $\{\hat{L}_1, \ldots, \hat{L}_k \}$. 

Then each end is isolated and so corresponds to a small open annulus on the punctured surface $S'$. 
Note that the surface $M'/\hat{\mathcal{F}} ' = S'$ can be identified with a subset of the class space $M/\hat{\mathcal{F}}$. 
This means that for any class $\hat{L}_i$, there is a small open disk $D_i$ in $M/\hat{\mathcal{F}}$ centered at $\hat{L}_i$ such that the intersection $D_i \cap S'$ is an open annulus. 
Thus the end compactification of $M'/\hat{\mathcal{F}}'$ is a surface with corners which is homeomorphic to the class space $M/\hat{\mathcal{F}}$. 
\end{proof}

Recall that each continuous codimension two compact foliation on a compact manifold is $R$-closed, and that the leaf (class) space of a continuous codimension two compact foliation of a compact manifold is a compact $2$-orbifold, which is homeomorphic to a surface with corners \cite{EMS,E,E2,V,V2}. 
In particular, the quotient space of a surface with corners by a finite group-action is also a surface with corners, which is homeomorphic to an orbifold as a topological space. 
Hence the following corollary is a generalization of the fact for codimension two compact foliations.  

\begin{corollary}\label{cor86}
The leaf class space of a codimension-two-like $R$-closed foliation on a compact manifold is a surface with corners. 
\end{corollary}

Corollary~\ref{cor86} implies (2) of Theorem~\ref{main:a}.  
Applying Lemma~\ref{lem:4.2} to the group-action cases, we have the following statement. 

\begin{lemma}\label{lem87}
Let $G$ be an $R$-closed group-action on a compact manifold $M$ and $H$ a finite index normal subgroup of $G$. 
If $H$ is codimension-two-like, then the orbit class space $M/\hat{\mathcal{F}} _G$ is a surface with corners. 
\end{lemma}

\begin{proof} 
By Theorem~\ref{main:c}, we have that $H$ is also $R$-closed. 
Applying Lemma~\ref{lem:4.2}, $M/\hat{\mathcal{F}}_H$ is a surface with corners. 
Since $G/H$ is a finite group acting $M/\hat{\mathcal{F}}_H$ and the quotient space of a surface by a finite group-action is a surface with corners, Lemma~\ref{lem:synd} implies that $M/\hat{\mathcal{F}} _G \cong (M/\hat{\mathcal{F}} _H)/(G/H)$ is a surface with corners. 
\end{proof}

Lemma~\ref{lem87} implies (2) of Theorem~\ref{main:b}. 

%

\section{another inherited properties of decompositions}

The double operation of a decomposition on a manifold preserves pointwise almost periodicity and $R$-closedness. 

\begin{lemma}\label{lem:decomp_double}
Let $\mathcal{F}$ be a decomposition of a manifold $M$ such that each boundary component consists of elements of $\mathcal{F}$, $D := M \cup - M$ the double of $M$, where $-M$ is the copy of $M$. 
Suppose that the boundary $\partial M$ consists of elements of $\mathcal{F}$. 
Define the induced decomposition $\mathcal{F}_D := \{ L, -L \mid L \in \mathcal{F} \}$ on $D$, called the double of $\mathcal{F}$, where $-L$ is a copy of $L$ in $-M$. 
Then the following statements hold: 
\\
$(1)$ The decomposition $\mathcal{F}$ is pointwise almost periodic if and only if so is the double $\mathcal{F}_D$. 
\\
$(2)$ The decomposition $\mathcal{F}$ is $R$-closed if and only if so is the double $\mathcal{F}_D$. 
\end{lemma}

\begin{proof}
Since each boundary component consists of elements of $\mathcal{F}$, any closures of element of $\mathcal{F}$ is contained in either the interior or the boundary of $M$. 
Then either $\hat{\mathcal{F}}(x) \subseteq \mathrm{int} M$ or $\hat{\mathcal{F}}(x) \subseteq \partial M$ for any point $x \in M$. 
Therefore either $\hat{\mathcal{F}}_D(x) \subseteq D - \partial M$ or $\hat{\mathcal{F}}_D(x) \subseteq \partial M$ for any point $x \in D$. 
This implies that $\hat{\mathcal{F}}(x) = \hat{\mathcal{F}}_D(x)$ for any point $x \in M$. 
This implies that the assertion $(1)$ holds. 
Since any subspace of a Hausdorff space is Hausdorff, the $R$-closedness of $\mathcal{F}_D$ implies one of $\mathcal{F}$. 
Taking the doubles of neighborhoods, the $R$-closedness of $\mathcal{F}$ implies one of $\mathcal{F}_D$. 
\end{proof}

\subsection{Semi-conjugacy for decompositions} 

Let $\mathcal{F}$ be a decomposition on a topological space $X$ and $\mathcal{G}$ a decomposition on a topological space $Y$. 
We call that a continuous surjection $f \colon X \to Y$ is a semi-conjugacy from $\mathcal{F}$ to $\mathcal{G}$ if the image of each element of $\F$ is an element of $\mathcal{G}$ (i.e. $f(\F(x)) = \mathcal{G}(f(x))$ for any $x \in X$).  
Then $\mathcal{F}$ is said to be semi-conjugate to $\mathcal{G}$. 
The semi-conjugacy preserves closures of elements. 

\begin{lemma}\label{lem:closure}
Let $\mathcal{F}$ be a  decomposition of a compact  space $X$, $\mathcal{G}$ a decomposition of a Hausdorff space $Y$,  and $p \colon  X \to Y$ a semi-conjugacy from $\mathcal{F}$ to $\mathcal{G}$. 
For any point $x \in X$, we have $p(\overline{\mathcal{F}(x)}) = \overline{\mathcal{G}(p(x))}$. 
\end{lemma}

\begin{proof}
Fix a point $x \in X$. 
Put $F := \mathcal{F}(x)$ and $L := \mathcal{G}(p(x))$. 
The semi-conjugacy implies that $p(F) = L$. 
Then $\overline{F} \subseteq p^{-1}(\overline{L})$ and so $L \subseteq p(\overline{F}) \subseteq \overline{L}$. 
Since $X$ is compact, so is $\overline{F}$. 
The Hausdorff property of $Y$ implies that  the compact subset $p(\overline{F})$ is closed  and so $\overline{L} = p(\overline{F})$.  
This means that $p(\overline{\mathcal{F}(x)}) = p(\overline{F}) = \overline{L} = \overline{\mathcal{G}(p(x))}$. 
\end{proof}

The semi-conjugacy preserves pointwise almost periodicity. 

\begin{lemma}\label{lem040}
Let $\mathcal{F}$ be a  decomposition of a compact  space $X$, $\mathcal{G}$ a decomposition of a Hausdorff space $Y$,  and $p \colon  X \to Y$ a semi-conjugacy from $\mathcal{F}$ to $\mathcal{G}$. 
If $\mathcal{F}$ is pointwise almost periodic, then so is $\mathcal{G}$. 
\end{lemma}

\begin{proof}
Assume that there is an element $L \in \mathcal{G}$ whose closure is not minimal. 
Then there is an element $L_0 \in \mathcal{G}$ such that $\overline{L_0} \subsetneq \overline{L}$. 
By the surjectivity of $p$, there is an element $F \in \mathcal{F}$ such that $p(F) = L$. 
Lemma~\ref{lem:closure} implies that $L_0 \subseteq \overline{L_0} \subsetneq \overline{L} = p(\overline{F})$ and so that $p^{-1}(L_0) \cap \overline{F} \neq \emptyset$. 
Fix a point $x_0 \in p^{-1}(L_0) \cap \overline{F}$ and put $F_0 := \mathcal{F}(x_0)$. 
Since $\F$ is pointwise almost periodic, the closure of any element of  $\F$ is invariant. 
Therefore $F_0 \subset \overline{F}$. 
Since $p(x_0) \in L_0$, the semi-conjugacy $p$ implies that $p(F_0) = p(\mathcal{F}(x_0)) = \mathcal{G}(p(x_0)) = L_0$. 
Lemma~\ref{lem:closure} implies that $p(\overline{F_0}) = \overline{L_0}  \subsetneq \overline{L} = p(\overline{F})$ and so $\overline{\mathcal{F}(x_0)}= \overline{F_0} \subsetneq \overline{F}$, which contradicts to the pointwise almost periodicity of $\mathcal{F}$. 
Thus $\mathcal{G}$ is pointwise almost periodic. 
\end{proof}

We show that a semi-conjugacy from $\mathcal{F}$ to $\mathcal{G}$ implies the a semi-conjugacy from $\hat{\mathcal{F}}$ to $\hat{\mathcal{G}}$ if the decompositions $\mathcal{F}$ and $\mathcal{G}$ are pointwise almost periodic. 

\begin{lemma}\label{lem040a}
Let $\mathcal{F}$ be a pointwise almost periodic decomposition of a compact  space $X$, $\mathcal{G}$ a decomposition of a Hausdorff space $Y$,  and $p \colon  X \to Y$ a semi-conjugacy from $\mathcal{F}$ to $\mathcal{G}$. 
Then $p(\hat{\F}(A)) = \hat{\mathcal{G}}(p(A))$ for any $A \subseteq X$. 
\end{lemma}

\begin{proof}
Lemma~\ref{lem040} implies that $\mathcal{G}$ is pointwise almost periodic. 
Since $\hat{\F}(A) = \bigcup_{x \in A} \hat{\F}(x)$, it suffices to show $p(\hat{\F}(x)) = \hat{\mathcal{G}}(p(x))$ for any $x \in X$. 
Fix any point $x \in X$. 
Lemma~\ref{lem:closure} implies that $p(\overline{\F(x)}) = \overline{\G(p(x))}$. 
%
The pointwise almost periodicity implies $p(\hat{\F}(x)) = p(\overline{\F(x)}) = \overline{\G(p(x))} = \hat{\G}(p(x))$.
\end{proof}

The semi-conjugacy preserves $R$-closedness.

\begin{proposition}\label{lem041}
Let $\mathcal{F}$ $(\mathrm{resp.}$ $\mathcal{G})$ be a decomposition of a compact Hausdorff space $X$ $(\mathrm{resp.}$ $Y)$. 
Suppose $\mathcal{F}$ is semi-conjugate to $\mathcal{G}$. 
If $\mathcal{F}$ is $R$-closed, then so is $\mathcal{G}$. 
\end{proposition}

\begin{proof}
Fix any closed subset $K$ of $Y$. 
By Theorem~\ref{th:usc}, it suffices to show that the saturation $\hat{\mathcal{G}}(K) = \mathcal{G}(K)$ of $\mathcal{G}$ is closed. 
Let $p \colon  X \to Y$ be a semi-conjugacy. 
Since $p^{-1}(K)$ is closed, the $R$-closedness of $\F$ implies that $\hat{\mathcal{F}}(p^{-1}(K))$ is closed and so compact, because $X$ is compact. By the surjectivity of $p$, we have  $p(p^{-1}(K)) = K$. 
%
The semi-conjugacy implies that $\hat{\mathcal{G}}(K)= \hat{\mathcal{G}}(p(p^{-1}(K))) = p(\hat{\mathcal{F}}(p^{-1}(K)))$ is compact and so closed, because $Y$ is Hausdorff. 
\end{proof}

\subsection{Semi-conjugacy and lifting properties for group actions} 
Consider semi-conjugacies for group-actions to apply Lemma~\ref{lem040} and Proposition~\ref{lem041}. 
For any group-action $G$ (resp. $H$) on a topological space $X$ (resp. $Y$), we call that a surjection $p \colon X \to Y$ is a semi-conjugacy if the image of an orbit of $G$ is an orbit of $H$ (i.e. $p(G(x)) = H(p(x))$ for any point $x \in X$). 
Then $G$ is said to be semi-conjugate to $H$. 
By Lemma~\ref{lem040} and Proposition~\ref{lem041}, we have the following statement. 

\begin{corollary}\label{cor042}
Suppose a group-action $G$ on a compact Hausdorff space $X$ is semi-conjugate to a group-action $H$ on a compact Hausdorff space $Y$. 
If $G$ is pointwise almost periodic $(\mathrm{resp.}$ $R$-closed$)$, then $H$ is pointwise almost periodic $(\mathrm{resp.}$ $R$-closed$)$. 
\end{corollary}

Notice that any flow generated by a vector field $X$ on a non-orientable connected compact manifold induces the induced flow on the orientable double cover because $X$ can be lifted to the cover.  
We also state the inherited property for orientation double covers to characterize $R$-closedness for homeomorphisms isotopic to the identity on non-orientable surfaces in Section~\ref{sec:nonori}.

\begin{proposition}\label{lem43+}
Let $G$ be a group-action on a non-orientable compact manifold $X$, $\widetilde{X}$ the orientable double cover of $X$. 
Suppose that there is a lift of $G$ on $\widetilde{X}$. 
Then $G$ is $R$-closed if and only if so is the lift of $G$ on $\widetilde{X}$. 
\end{proposition}

\begin{proof}  
Let $\F_G$ be the decomposition induced by $G$ and $p \colon  \widetilde{X} \to X$ the orientable double covering map. 
Denote by $\F_{\widetilde{G}}$ the decomposition induced by the lift $\widetilde{G}$ of $G$ on $\widetilde{X}$. 
Let $\Gamma(p)$ be the covering transformation group of the orientation double covering map $p$ and $\F_{\widetilde{G} \times \Gamma(p)}$ the decomposition induced by the product group $\widetilde{G} \times \Gamma(p)$. 
Then $\Gamma(p)$ is isomorphic to $\Z/2\Z$ and the class space $X/\hat{\F}_G$ can be identified with $(\widetilde{X}/\hat{\F}_{\widetilde{G}})/\Gamma(p) = \widetilde{X}/\hat{\F}_{\widetilde{G} \times \Gamma(p)}$. 
Moreover, we identify $\widetilde{G}$ with the finite index normal subgroup of $\widetilde{G} \times \Gamma(p)$. 

Suppose that the lift $\widetilde{G}$ of $G$ on $\widetilde{X}$ is $R$-closed. 
Applying Theorem~\ref{main:c}, the product group $\widetilde{G} \times \Gamma(p)$ is $R$-closed. 
Theorem~\ref{th:usc} implies that $\widetilde{X}/\hat{\F}_{\widetilde{G} \times \Gamma(p)} = X/\hat{\F}_G$ is Hausdorff and so $G$ is $R$-closed. 
Conversely, suppose that $G$ is $R$-closed. 
Theorem~\ref{th:usc} implies that $\widetilde{X}/\hat{\F}_{\widetilde{G} \times \Gamma(p)} = X/\hat{\F}_G$ is Hausdorff and so $\widetilde{G} \times \Gamma(p)$ is $R$-closed. 
Applying Theorem~\ref{main:c}, the finite index normal subgroup $\widetilde{G}$ is $R$-closed. 
\end{proof}

\section{Complete characterizations of $R$-closedness for homeomorphisms on orientable compact surfaces}\label{sec:ori}
In this section, we characterize $R$-closedness for homeomorphisms on orientable compact surfaces.

\subsection{Pointwise almost periodicity and $R$-closedness for homeomorphisms on compact surfaces}

A homeomorphism $f \colon  X \to X$ is periodic if there is a positive integer $n$ such that $f^n$ is identical. 
A homeomorphisms $f \colon S \to S$ is semi-conjugate to a homeomorphism $g \colon  T \to T$ if there is a surjection $p \colon S \to T$ with $p \circ f = g \circ p$. 
Then $p$ is called a semi-conjugacy. 
Note that if a homeomorphisms $f \colon S \to S$ is semi-conjugate to a homeomorphism $g \colon  T \to T$ via a semi-conjugacy $p \colon S \to T$, then the image of an orbit of $f$ is an orbit of $g$. 
A minimal set $\mathcal{C}$ on a surface homeomorphism $f  \colon  S \to S$ is an extension of a Cantor set (resp. a periodic orbit) if there are a homeomorphism $g \colon  S \to S$ and a semi-conjugacy $p \colon  S \to S$ with $p \circ f = g \circ p$ such that $p$ is homotopic to the identity  and $p(\mathcal{C})$ is a Cantor set (resp. a periodic orbit) which is a minimal set of $g$. 
We show the following non-existence of a minimal set which is an extension of a Cantor set for an $R$-closed homeomorphism on an orientable compact surface.

\begin{lemma}\label{lem55}
Let $f \colon  S \to S$ be an $R$-closed homeomorphism on an orientable compact surface $S$.   
Then $f$ has no minimal set which is an extension of a Cantor set.
\end{lemma}

\begin{proof}
Assume that there is a minimal set $\mathcal{M}$ of $f$ which is an extension of a Cantor set.
By Lemma~\ref{lem:decomp_double}, taking the double of $S$, we may assume that $S$ is closed. 
By the definition of an extension of a Cantor set, there are a surface homeomorphism $g \colon  S \to S$, a semi-conjugacy $p \colon  S \to S$ from $f$ to $g$, and a Cantor minimal set $\mathcal{C} := p(\mathcal{M})$ of $g$. It suffices to show the non-existence of such a minimal set $\mathcal{C}$. 
By Corollary~\ref{cor042}, we obtain that $g$ is an $R$-closed homeomorphism with a Cantor minimal set $\mathcal{C}$ on a closed orientable surface $S$. 
Since a Cantor set is zero-dimensional, there is a basis of $\mathcal{C}$ consisting of closed and open subsets. 

We claim that, for a point $x \in \mathcal{C}$, there is an open disk $U_x$ which is a neighborhood of $x$ and whose boundary $\partial U_x$ is a circle but does not intersect $\mathcal{C}$. 
Indeed, fix a small open ball $U_0$ centered at $x$. 
Since there is a basis of $\mathcal{C}$ consisting of closed and open subsets, there is an open \nbd $U_1 \subseteq U_0$ such that $\mathcal{C} \cap \partial U_1 = \emptyset$. 
Recall that the filling $\mathrm{Fill}_{U_0}(U_1)$ of $U_1$ in $U_0$ is defined as follows: $y \in \mathrm{Fill}_{U_0}(U_1)$ if either $y \in B$ or there is an open disk in $U_0$ containing $x$ whose boundary is contained in $U_1$. 
Then the filling $U_x := \mathrm{Fill}_{U_0}(U_1)$ is a desired open disk containing $x$ such that $\mathcal{C} \cap \partial U_x = \emptyset$.  

The compactness of $\mathcal{C}$ implies that there are finitely many points $x_1, \ldots, x_k \in \mathcal{C}$ with $\mathcal{C} \subset \bigcup_{i=1}^k U_{x_i}$.  
Let $D_{x_i} := \F_g(\partial U_{x_i})$ be the saturation of the boundary of $U_{x_i}$ and $D := \bigcup_{i=1}^k D_{x_i}$. 
Then $D \cap \mathcal{C} = \emptyset$. 
By the $R$-closedness of $g$, the saturation $D$ is closed. 
Since $X$ is normal, there are an disjoint open neighborhoods $U_2 \subseteq \bigcup_{i=1}^k U_{x_i} \cap g(\bigcup_{i=1}^k U_{x_i})$ and $V_2$ of $\mathcal{C}$ and $D$ respectively with $V_2 = \bigcup_{i=1}^k V_{x_i}$ such that $V_{x_i}$ is an annular neighborhood of $\partial U_{x_i}$ and each $U_{x_i} \setminus V_{x_i}$ is a closed disk.  
Notice that each connected component of $U_2$ is contained in some closed disk $U_{x_i} \setminus V_{x_i}$ and so contractible in $U_{x_i} \setminus V_{x_i}$.  
By the $R$-closedness of $g$, there is a $g$-invariant open neighborhood $U \subseteq U_2$ of $\mathcal{C}$. 
Then each connected component of $U$ is contained and contractible in some closed disk $U_{x_i} \setminus V_{x_i}$. 
Fix any point $x \in \mathcal{C}$. 
Let $W$ be the connected component of $U$ containing $x$ and $B_x := U_{x_l} \setminus V_{x_l}$ the closed disk containing $W$ for some $l =1, \ldots, k$. 
Then the connected component $W$ is an open connected subset in the closed disk $B_x$. 
Therefore the filling $\mathrm{Fill}_{B_x}(W) \subset B_x$ of $W$ is an open disk. 
Since $\mathcal{C} \cap D = \emptyset$, the $g$-invariance of $U$ implies that the connected component $W$ of $U$ is the connected component of $B_x \cap \bigcup_{n \in \mathbb{Z}} g^n(W)$ containing $x$. 
Then $\partial W \subseteq \partial \F_g(W)$. 
Define the first return map $G  \colon W \to W$ by $G(y) := g^{N(y)}(y)$ for any $y \in W$, where ${N(y)} := \min \{ n \in \mathbb{Z}_{>0} \mid g^n(y) \in W \}$. 
The connectivity of $W$ implies that $N(y)$ is a constant $N(W) > 0$, and that $g^{\pm N(W)} (W) \subseteq W$.  
This means that $G(W) = g^{N(W)} (W) = W$ and so that $G$ is a homeomorphism. 
Note that $x \in W \subset B_x$.  
Since $U$ is $g$-invariant, for any open disk $D_\gamma$ in $\mathrm{Fill}_{B_x}(W) \subseteq B_x$ with $\gamma := \partial D_\gamma \subseteq \overline{W}$ and for any $n \in \mathbb{Z}$, the image $G^n(\gamma) \subseteq \mathrm{Fill}_{B_x}(W)$ is the boundary of an open disk in $\mathrm{Fill}_{B_x}(W) \subseteq B_x$. 
By $\overline{W} \subseteq g^{\pm N(W)}(\overline{W})$, we have $\overline{W} \subseteq g^{\mp N(W)}(\overline{W})$ and so $\overline{W} = g^{N(W)}(\overline{W})$. 
Therefore we can extend $G$ into a homeomorphism $G := g^{N(W)} \colon \overline{W} \to \overline{W}$ because $g^{N(W)}|_{\overline{W}}$ is bijective and $\overline{W}$ is compact Hausdorff. 

We claim that $G$ maps the outer boundary $\partial (\mathrm{Fill}_{B_x}(W)) \subseteq \partial W$ of $W$ into itself.  
Indeed, assume that the outer boundary $\partial_{\mathrm{out}} := \partial (\mathrm{Fill}_{B_x}(W)) \subseteq \partial W$ is not preserved by $G$. 
Since $G  \colon W \to W$ is homeomorphic, there is a connected component $\partial_{\mathrm{in}} \neq \partial_{\mathrm{out}}$ of $\partial W$ such that $G(\partial_{\mathrm{in}}) = \partial_{\mathrm{out}}$. 
Then the surface $S$ is covered by two open disks $\mathrm{int} (\mathrm{Fill}_{B_x}(W))$ and $G(\mathrm{int} (\mathrm{Fill}_{B_x}(W)))$ such that the intersection $\mathrm{int} (\mathrm{Fill}_{B_x}(W)) \cap G(\mathrm{int} (\mathrm{Fill}_{B_x}(W)))$ is an open annulus. 
This means that $S$ is a sphere and so there is a spherical $R$-closed homeomorphism with Cantor minimal sets, which contradicts to \cite[Corollary~2.6]{Y2}. 

Then $\partial (\mathrm{Fill}_{B_x}(W))$ is an invariant closed subset with respect to $G$. 
By the construction of $\mathrm{Fill}_{B_x}(W)$, the injectivity of $G$ implies that we can extend $G$ into a mapping $G  \colon \overline{\mathrm{Fill}_{B_x}(W)} \to \overline{\mathrm{Fill}_{B_x}(W)}$. 
Since the original $g$ is $R$-closed and homeomorphic, by Theorem~\ref{main:c}, so is the iteration $G$. 
By Theorem~\ref{th:usc}, since $G$ is the first return map of $W$, the quotient space $\mathrm{Fill}_{B_x}(W)/\hat{G} \cong \F_g(\mathrm{Fill}_{B_x}(W))/\hat{g}$ is Hausdorff, where the relation $X_1 \cong X_2$ denotes that a topological space $X_1$ is homeomorphic to a topological space $X_2$. 
%
%
By the one-point compactification $S_\infty$ of $\mathrm{Fill}_{B_x}(W)$, the resulting surface $S$ is a sphere. 
By the compactness of $\partial (\mathrm{Fill}_{B_x}(W))$, the $R$-closedness of $G$ and the normality of $S$ imply that the boundary $\partial (\mathrm{Fill}_{B_x}(W))$ and the orbit class $\hat{G}(b)$ for any point $b \in \mathrm{Fill}_{B_x}(W)$ can be separated by disjoint invariant open neighborhoods. 
Adding a new fixed point, we obtain the resulting homeomorphism $G_\infty$ on $S_\infty$ induced by $G$.
This implies that $S_\infty$ is homeomorphic to $\overline{\mathrm{Fill}_{B_x}(W)}/ \sim_{\partial}$ and that the orbit class space of the induced homeomorphism $G_\infty$ on the sphere $S_\infty$ is Hausdorff, where the equivalent relation $\sim_{\partial}$ is defined as follows: $ x \sim_{\partial} y $ if either $x = y$ or $\{x, y \} \subseteq \partial (\mathrm{Fill}_{B_x}(W))$.    
By Theorem~\ref{th:usc}, the resulting homeomorphism $G_\infty$ is an $R$-closed homeomorphism on a sphere with a Cantor minimal set. 
This contradicts to \cite[Corollary~2.6]{Y2}. 
Thus there are no Cantor minimal sets. 
\end{proof}

We show the following statement. 

\begin{lemma}\label{disjoint}
The orbit class space of a pointwise almost periodic toral homeomorphism each of whose minimal set is a finite disjoint union of essential circloids is either a circle or a closed interval. 
In particular, such a homeomorphism is $R$-closed. 
\end{lemma}

\begin{proof}
Let $f \colon  \T^2 \to \T^2$ be a pointwise almost periodic toral homeomorphism each of whose minimal set is a finite disjoint union of essential circloids. 
Fix such an essential circloid $\gamma$. 
Since circloids are annular continua, there is an open annular \nbd $U$ such that $U - \gamma$ is a disjoint union of two open essential annuli. 

We claim that the complement $\T^2 - \gamma$ is an open annulus such that any essential simple closed curves in the union $\T^2 - \gamma = (\T^2 - U) \cup (\overline{U} - \gamma)$ are parallel to each other. 
Indeed, replacing an open annular \nbd of $\gamma$ contained in $U$, we may assume that the boundary $\partial U$ consists of disjoint two simple closed curves. 
Then the complement $\T^2 - \overline{U}$ is an open annulus and so any essential simple closed curves in the union $\T^2 - \gamma = (\T^2 - U) \cup (\overline{U} - \gamma)$ are parallel to each other. 
This means that the submanifold $\T^2 - \gamma = (\T^2 - U) \cup (\overline{U} - \gamma)$ is an open annulus. 

Therefore the set of connected components of minimal sets contained in $U$ is totally ordered. 
Collapsing connected components of minimal sets into singletons, the resulting space $\T^2/\hat{f}$ is a circle and so the induced mapping on the circle is a homeomorphism whose orbits are finite and so pointwise periodic. 
By \cite[Theorem]{M}, each pointwise periodic homeomorphism on a circle is periodic and so is the induced mapping and so the orbit class space of $f$ is a closed interval if the induced mapping is non-orientable, and it is a circle if the induced mapping is orientable. 
Theorem~\ref{th:usc} implies that $f$ is $R$-closed. 
\end{proof}

\subsection{Characterize $R$-closedness for homeomorphisms on orientable surfaces}

It is known that a non-periodic spherical homeomorphism is $R$-closed if and only if the orbit class space is a closed interval (see \cite[Theorem~6 and Theorem~7]{Ma}, and \cite[Corollary~2.6]{Y2}). 
Moreover, we state an another characterization as follows. 

\begin{proposition}\label{sphere}
An orientation-preserving {\rm(resp.} orientation-reversing{\rm)} pointwise almost periodic homeomorphism $f$ on a sphere is $R$-closed if and only if either $f$ is periodic or the sphere consists of two fixed points and minimal sets of $f$ {\rm(resp.} $f^2${\rm)} each of which is a circloid. 
\end{proposition}

\begin{proof}
Suppose that $f$ is orientation-preserving. 
Since the orbit class space of a periodic homeomorphism on a compact manifold is Hausdorff, we may assume that $f$ is not periodic. 
Suppose that $f$ is $R$-closed. 
 \cite[Corollary~2.6]{Y2} implies that the sphere consists of two fixed points and minimal sets each of which is a circloid. 
Conversely, suppose that the sphere consists of two fixed points and minimal sets each of which is a circloid. 
Fix a circloid $\gamma$. 
Since circloids are annular continua, there is an open annular \nbd $U$ such that $U - \gamma$ is a disjoint union of two open annuli. 

We claim that the complement $\mathbb{S}^2 - \gamma$ is the disjoint union of two open disks. 
Indeed, replacing an open annular \nbd of $\gamma$ contained in $U$, we may assume that the boundary $\partial U$ consists of disjoint two simple closed curves. 
Denote by $\gamma_i$ the boundary components of $U$ ($i=1,2$). 
Recall that Jordan-Schoenflies theorem says that the complement of any simple closed curve in a sphere is a disjoint union of two open disks. 
Therefore the complement $\mathbb{S}^2 - \gamma_1$ is a disjoint union of two open disks. 

Let $D$ be the open disk containing $U \sqcup \gamma_2$. 
Jordan curve theorem implies that the complement $D - \gamma_2$ is a disjoint union of an open disk and an open annulus. 
This means that the complement $\mathbb{S}^2 - (\gamma_1 \sqcup \gamma_2) = \mathbb{S}^2 - \partial U$ consists two open disks $D_1$, $D_2$ and one open annulus $U$ such that $\partial D_i = \gamma_i$ and $\partial U = \gamma_1 \sqcup \gamma_2$. 
Let $U_i$ be the connected components of the complement $U - \gamma_i$ such that one of boundary components is $\gamma_i$ ($i=1,2$). 
Then each union $S_i := D_i \sqcup \gamma_i \sqcup U_i$ is open and simply connected. 
Since the Riemann mapping theorem states that each nonempty open simply connected proper subset of $\mathbb{C}$ is conformally equivalent to the unit disk, the unions $S_1$ and $S_2$ are open disks such that $\mathbb{S}^2 = S_1 \sqcup \gamma \sqcup S_2$. 
Thus the orbit class space is a closed interval and so $f$ is $R$-closed. 

Suppose that $f$ is orientation-reversing. 
Then $f^2$ is orientation-preserving. 
By Theorem~\ref{main:c}, we have that $f$ is $R$-closed if and only if so is $f^2$. 
Therefore the characterization of $R$-closedness for orientable spherical homeomorphisms implies the assertion. 
\end{proof}

\subsection{Proofs of Theorem~\ref{main:d}  and Theorem~\ref{main:e}}
\subsubsection{Proof of Theorem~\ref{main:d}}
%
%
Suppose that the homeomorphism $f$ is $R$-closed. 
Assume that $f$ is neither minimal nor periodic. 
By Lemma~\ref{lem55} and  \cite[Theorem~3.2]{Y2}, each minimal set is a finite disjoint union of essential circloids. 
Conversely, suppose that $f$ is minimal.
Then the orbit class space is a singleton and so is Hausdorff. 
Theorem~\ref{th:usc} implies that $f$ is $R$-closed. 
Suppose that $f$ is periodic. 
Since the leaf space of a continuous codimension two compact foliation of a compact manifold is a compact orbifold \cite{EMS,E,E2,V,V2} and so is Hausdorff, Theorem~\ref{th:usc} implies that $f$ is $R$-closed. 
Suppose that $f$ is pointwise almost periodic and each minimal set is a finite disjoint union of essential circloids. 
Lemma~\ref{disjoint} implies that $f$ is $R$-closed. 

\subsubsection{Proof of Theorem~\ref{main:e}}
%
%
Let $f$ be a homeomorphism on an orientable compact surface $S$ whose Euler number is negative. 
Suppose that $f$ is periodic. 
Since the orbit class space of a periodic homeomorphism on a compact manifold is Hausdorff, the homeomorphism $f$ is $R$-closed. 
Conversely, suppose that $f$ is $R$-closed. 
By Lemma~\ref{lem:decomp_double}, taking the double of $S$, we may assume that $S$ is closed and the Euler number of $S$ is negative. 
Then the genus of $S$ is at least two.   
By Lemma~\ref{lem55}, there is no minimal set which is an extension of a Cantor set.
By  \cite[Theorem~3.1]{Y2}, we have that $f$ is pointwise periodic. 
By Theorem~\cite{M}, each pointwise periodic surface homeomorphism is periodic and so is $f$. 

\subsection{Characterize $R$-closedness for homeomorphisms on closed disks and annuli}

We characterize $R$-closedness for homeomorphisms on closed disks. 

\begin{proposition}\label{prop:611}
Each homeomorphism on a closed disk is $R$-closed if and only if  it satisfies one of the following statements:
\\
$(1)$ The homeomorphism $f$ is periodic. 
\\
$(2)$ The homeomorphism $f$ is orientation-preserving and the disk consists of one fixed point and minimal sets each of which is a circloid. 
\end{proposition}

\begin{proof}
Let $f \colon \mathbb{D} \to \mathbb{D}$ be a homeomorphism on a closed disk $\mathbb{D}$. 
By $f(\partial \mathbb{D}) = \partial \mathbb{D}$, the orientation-reversing property of $f$ implies that the restriction $f|_{\partial \mathbb{D}}$ is a homeomorphism on a circle $\partial \mathbb{D}$. 
Let $\mathbb{S}^2$ be the double of $\mathbb{D}$. 
Since the boundary $\partial \mathbb{D}$ is invariant under $f$, the induced mapping $f_{\mathbb{S}^2} \colon \mathbb{S}^2 \to \mathbb{S}^2$ is well-defined and homeomorphic such that $f_{\mathbb{S}^2}|_\mathbb{D} = f$. 
Then $f$ is periodic if and only if so is $f_{\mathbb{S}^2}$. 
Therefore Lemma~\ref{lem:decomp_double} implies that the periodicity of $f$ implies the $R$-closedness of $f$. 
Thus we may assume that $f$ is not periodic. 
Suppose that $f$ is orientation-preserving and the disk consists of one fixed point and minimal sets each of which is a circloid. 
The invariance $\mathbb{D}$ of $f$ implies that the double $\mathbb{S}^2$ consists of two fixed points and minimal sets each of which is a circloid. 
By Proposition~\ref{sphere}, the induced mapping $f_{\mathbb{S}^2}$ is an $R$-closed homeomorphism. 
Lemma~\ref{lem:decomp_double} implies that $f$ is $R$-closed. 

Conversely, suppose that $f$ is $R$-closed. 
Lemma~\ref{lem:decomp_double} implies that $f_{\mathbb{S}^2}$ is $R$-closed. 
We claim that $f$ is orientation-preserving. 
Indeed, assume that $f$ is orientation-reversing. 
Then so are $f_{\mathbb{S}^2}$ and $f|_{\partial \mathbb{D}}$. 
The non-orientability of $f|_{\partial \mathbb{D}}$ implies fixed points for the homeomorphism $f|_{\partial \mathbb{D}}$.   
Proposition~\ref{sphere} implies that the sphere $\mathbb{S}^2$ consists of two fixed points and minimal sets of $f_{\mathbb{S}^2}$ each of which is a circloid. 
By the invariance of $\mathbb{D}$, the disk $\mathbb{D}$ consists of a fixed point and minimal sets of $f^2$ each of which is a circloid. 
The invariance of $\partial \mathbb{D}$ implies that the restriction $f^2|_{\partial \mathbb{D}}$ is minimal, which contradicts the existence of fixed points of $f|_\mathbb{\partial D}$. 

By the orientability of $f$, Proposition~\ref{sphere} implies that the sphere $\mathbb{S}^2$ consists of two fixed points and minimal sets of $f_{\mathbb{S}^2}$ each of which is a circloid. 
The invariance of $\partial \mathbb{D}$ implies that $f(\mathbb{D}) = \mathbb{D}$ and $f(\partial \mathbb{D}) = \partial \mathbb{D}$. 
This means that the disk $\mathbb{D}$ consists of a fixed point and minimal sets of $f$ each of which is a circloid. 
\end{proof}

We characterize $R$-closedness for homeomorphisms on closed annuli. 

\begin{proposition}\label{prop:annulus}
An orientation-preserving {\rm(resp.} orientation-reversing{\rm)} homeomorphism on a closed annulus is $R$-closed if and only if it satisfies one of the following statements:
\\
$(1)$ The homeomorphism $f$ is periodic. 
\\
$(2)$ The annulus consists of minimal sets of $f$ {\rm(resp.} $f^2${\rm)} each of which is a circloid. 
\end{proposition}

\begin{proof}
Let $f \colon \mathbb{A} \to \mathbb{A}$ be a homeomorphism on a closed annulus $\mathbb{A}$. 
Then $f(\partial \mathbb{A}) = \partial \mathbb{A}$. 
Since the orbit class space of a periodic homeomorphism on a compact manifold is Hausdorff, if $f$ is periodic then $f$ is $R$-closed. 
Thus we may assume that $f$ is not periodic. 
Denote by $\partial_1$ and $\partial_2$  the boundary components of $\mathbb{A}$. 
Then $f^2(\partial_i) = \partial_i$ for any $i = 1,2$. 
Denote by $D_1$ and $D_2$ closed disks. 
Attaching $\partial_i$ to $\partial D_i$, the resulting surface $\mathbb{S}^2 := \mathbb{A} \cup (D_1 \sqcup D_2)$ is a sphere. 
Let $S$ be the resulting space from $\mathbb{S}^2$ collapsing disks $D_i$ into singletons. 

We claim that the induced mapping $f_{S} \colon S \to S$ by $f$ is well-defined and is a homeomorphism. 
Indeed, Moore's theorem (cf. p.3 in \cite{D}) implies that the quotient space $S$ is a sphere. 
By construction, the sphere $S$ can be identified with the resulting space from the closed annulus $\mathbb{A}$ collapsing boundaries $\partial_i = \partial D_i$ into singleton. 
The invariance of $\partial_i$ implies that the resulting mappings $f_S \colon S \to S$ (resp. $f_S^{-1}$) from $f \colon \mathbb{A} \to \mathbb{A}$ (resp. $f^{-1}$) is well-defined and homeomorphic. 

By construction, the orbit class space $\mathbb{A}/f$ is homeomorphic to the orbit class space $S/f_S$. 
Suppose that $f$ is $R$-closed. 
Then so is $f_S$. 
By Proposition~\ref{sphere}, if $f$ is orientation-preserving (resp. orientation-reversing) then the sphere $S$ consists of two fixed points and minimal sets of $f_S$ (resp. $f_S^2$) each of which is a circloid. 
This means that $\mathbb{A}$ consists of minimal sets of $f$ (resp. $f^2$) each of which is a circloid. 
Conversely, putting $g:= f$ (resp. $g:=f^2$) if $f$ is orientation-preserving (resp. orientation-reversing), suppose that $\mathbb{A}$ consists of minimal sets of $g$ which are circloids. 
Then $g$ is orientation-preserving such that each boundary component $\partial_i$ is a minimal set of $g$. 
Therefore the sphere $S$ induced from $\mathbb{A}$ consists of two fixed points and minimal sets of $g$ each of which is a circloid. 
By Proposition~\ref{sphere}, the induced mapping $g$ is $R$-closed. 
Theorem~\ref{main:c} implies that $f$ is $R$-closed if and only if so is $f^2$. 
This means that $f$ is $R$-closed.  
\end{proof}

\section{Complete characterizations of $R$-closedness for homeomorphisms isotopic to the identity on non-orientable compact surfaces}\label{sec:nonori}
In this section, we characterize $R$-closedness for homeomorphisms isotopic to the identity on non-orientable compact surface.

\subsection{Existence of lifts of homeomorphisms}
We state a necessary and sufficient condition for the existence of lifts of homeomorphisms. 

\begin{lemma}\label{lifting}
Let $p \colon  \widetilde{X} \to X$ be the covering map of a path-connected and locally path-connected space $X$, $h \colon  X \to X$ a homeomorphism, $\widetilde{x} \in \widetilde{X}$, and $\widetilde{y} \in \widetilde{X}$ with $p(\widetilde{y})=h \circ p (\widetilde{x})$. 
Then the following are equivalent:
\\
$(1)$ $h_*(p_*(\pi_1(\widetilde{X}, \widetilde{x}))) = p_*(\pi_1(\widetilde{X}, \widetilde{y}))$. 
\\
$(2)$ There is a lift $\widetilde{h} \colon  \widetilde{X} \to \widetilde{X}$ of $h$ with $\widetilde{h}(\widetilde{x}) = \widetilde{y}$. 
\\
Here $h_* \colon  \pi_1(X,  p(\widetilde{x})) \to \pi_1(X,h \circ p(\widetilde{x}))$ and $p_* \colon  \pi_1(\widetilde{X}, \widetilde{x}) \to \pi_1(X, p(\widetilde{x}))$ are the induced maps of $h$ and $p$ on the fundamental groups  respectively. 
In any case,  the lift $\widetilde{h} \colon  \widetilde{X} \to \widetilde{X}$ is a homeomorphism. 
\end{lemma}

\begin{proof}
Since $h$ is homeomorphic, the composition $f \colon = h \circ p \colon  \widetilde{X} \to X$ is a covering map and the induced map $h_* \colon  \pi_1(X,  p(\widetilde{x})) \to \pi_1(X,h \circ p(\widetilde{x}))$ is isomorphic. 
Put $x := p(\widetilde{x})$ and $y := h(x) = h \circ p (\widetilde{x}) = p(\widetilde{y})$. 

Suppose that there is a lift $\widetilde{h} \colon  \widetilde{X} \to \widetilde{X}$ of the homeomorphism $h$ with $\widetilde{h}(\widetilde{x}) = \widetilde{y}$. 
Then $p \circ \widetilde{h} = h \circ p$. 
The unique lifting property (cf.  \cite[Proposition~1.34]{Hat}) implies that the lift preserving a base point of an identity is identical. 
Therefore $\widetilde{h}$ is a homeomorphism because the compositions $\widetilde{h^{-1}} \circ \widetilde{h}$ and $\widetilde{h} \circ \widetilde{h^{-1}}$ are identical, where $\widetilde{h^{-1}} \colon  \widetilde{X} \to \widetilde{X}$ is the lift of the homeomorphism $h^{-1}$ with $\widetilde{h^{-1}}(\widetilde{y}) = \widetilde{x}$. 
This implies that $h_*(p_*(\pi_1(\widetilde{X}, \widetilde{x}))) = (h \circ p)_*(\pi_1(\widetilde{X}, \widetilde{x})) = (p \circ \widetilde{h})_*(\pi_1(\widetilde{X}, \widetilde{x}))) = p_*(  \widetilde{h}_*(\pi_1(\widetilde{X}, \widetilde{x})))) = p_*(\pi_1(\widetilde{X}, \widetilde{y}))$. 

Conversely, suppose that $h_*(p_*(\pi_1(\widetilde{X}, \widetilde{x}))) = p_*(\pi_1(\widetilde{X}, \widetilde{y}))$. 
Since $h$ is homeomorphism, a lifting criterion (cf.  \cite[Proposition~1.33]{Hat}) implies that there is a lift $\widetilde{h} \colon  \widetilde{X} \to \widetilde{X}$ of $h$ with $\widetilde{h}(\widetilde{x}) = \widetilde{y}$. 
Moreover, the restriction $h_*| \colon   p_*(\pi_1(\widetilde{X}, \widetilde{x})) \to  p_*(\pi_1(\widetilde{X}, \widetilde{y}))$ is isomorphic and so $h^{-1}_* (p_*(\pi_1(\widetilde{X}, \widetilde{y}))) = p_*(\pi_1(X, \widetilde{x}))$. 
This implies that there is a lift $\widetilde{h^{-1}} \colon  \widetilde{X} \to \widetilde{X}$ of $h^{-1}$ with $\widetilde{h^{-1}}(\widetilde{y}) = \widetilde{x}$. 
Then $p \circ \widetilde{h^{-1}} \circ \widetilde{h} = (h^{-1}) \circ p \circ \widetilde{h} = (h^{-1}) \circ h \circ p = p$ and so the composition $\widetilde{h^{-1}} \circ \widetilde{h} \colon  \widetilde{X} \to \widetilde{X}$ is a lift of the identity on $X$. 
The unique lifting property (cf.  \cite[Proposition~1.34]{Hat}) implies that $\widetilde{h}$ is a homeomorphism. 
\end{proof}

Recall the fact that homotopic homeomorphisms of surfaces are isotopic \cite{B,B2,E3}. 
We show that such homeomorphisms can be lifted to a covering spaces. 

\begin{corollary}\label{lifting_isotopic}
Let $p \colon  \widetilde{X} \to X$ be the covering map of a path-connected and locally path-connected space $X$ and $h \colon  X \to X$ a homeomorphism isotopic to the identity. 
Then there is a lift $\widetilde{h} \colon  \widetilde{X} \to \widetilde{X}$ of $h$ which is a homeomorphism isotopic to the identity. 
\end{corollary}

\begin{proof}
Denote by $H: [0,1] \times X \to X$ an isotopy from identity to $h$. 
Fix a point $\widetilde{x} \in \widetilde{X}$. 
Put $x := p(\widetilde{x})$ and $y := h(x)$. 
Fix a point $\widetilde{y} \in \widetilde{X}$ with $y = p(\widetilde{y})$. 

We claim that $h_*(p_*(\pi_1(\widetilde{X}, \widetilde{x}))) = p_*(\pi_1(\widetilde{X}, \widetilde{y}))$. 
Indeed, for any simple closed curve $\gamma$ representing an element $[\gamma]$ of the subgroup $p_*(\pi_1(\widetilde{X}, \widetilde{x}))$ of the fundamental group $\pi_1(X,x)$, by the lifting criterion to the restriction $H| \colon [0,1] \times \gamma \to X$, there is a lift $\widetilde{h(\gamma)}$ of $h(\gamma)$ and so $[h(\gamma)] = p_*([\widetilde{h(\gamma)}]) \in p_*(\pi_1(\widetilde{X}, \widetilde{y}))$. 
This means $h_*(p_*(\pi_1(\widetilde{X}, \widetilde{x}))) \leq p_*(\pi_1(\widetilde{X}, \widetilde{y}))$. 
Similarly, for any simple closed curve $\mu$ representing an element $[\mu]$ of the subgroup $p_*(\pi_1(\widetilde{X}, \widetilde{y}))$ of the fundamental group $\pi_1(X,y)$, by the lifting criterion to the restriction $H| \colon [0,1] \times h^{-1}(\mu) \to X$, there is a lift $\widetilde{h^{-1}(\mu)}$ of $h^{-1}(\mu)$ and so $[\mu] = h_*(p_*([\widetilde{h^{-1}(\mu)}]) \in h_*(p_*(\pi_1(\widetilde{X}, \widetilde{x})))$. 
This means $h_*(p_*(\pi_1(\widetilde{X}, \widetilde{x}))) \geq p_*(\pi_1(\widetilde{X}, \widetilde{y}))$. 

Since $H(1, \cdot ) = h$, this implies that $(H \circ (1_{[0,1]} \times p))_*(\pi_1([0,1] \times \widetilde{X}, \{1\} \times \widetilde{x}))) = H_*((1_{[0,1]} \times p)_*(\pi_1([0,1] \times \widetilde{X}, \{1\} \times \widetilde{x}))) = h_*(p_*(\pi_1(\widetilde{X}, \widetilde{x}))) = p_*(\pi_1(\widetilde{X}, \widetilde{y}))$, where $1_A \colon  A \to A$ is the identity on a topological space $A$ and the product map $1_{[0,1]} \times p$ is defined by $(1_{[0,1]} \times p)(t,x) = (t, p(x))$. 
 \cite[Proposition~1.33]{Hat} implies that there is a lift $\widetilde{H} \colon  [0,1] \times \widetilde{X} \to \widetilde{X}$ of the composition $H \circ (1_{[0,1]} \times p) \colon  [0,1] \times \widetilde{X} \to X$ with $\widetilde{H}(1, \widetilde{x}) = \widetilde{y}$. 
Then $p \circ \widetilde{H}(t, \cdot) = H (t, \cdot) \circ p$. 
Lemma~\ref{lifting} implies $\widetilde{H}(t, \cdot) \colon  \widetilde{X} \to \widetilde{X}$ for any $t \in [0,1]$ is a homeomorphism.  
This means that $\widetilde{H}$ is an isotopy with  $\widetilde{H}(0, \cdot) = 1_{\widetilde{X}}$ and $\widetilde{H}(1, \cdot ) = \widetilde{h}$. 
\end{proof}

\subsection{Characterize $R$-closedness for homeomorphisms isotopic to the identity on non-orientable surfaces}

We have the following characterization of $R$-closedness for a pointwise almost periodic but non periodic homeomorphism isotopic to the identity on a projective plane. 

\begin{proposition}
The following are equivalent for a pointwise almost periodic but non periodic homeomorphism $f$ isotopic to the identity on a projective plane $\mathbb{P}^2$: 
\\
$(1)$ The homeomorphism $f$ is $R$-closed. 
\\
$(2)$ The orbit class space is a closed interval. 
\\
$(3)$ The projective plane $\mathbb{P}^2$ consists of one fixed point and minimal sets each of which is lifted into a circloid by the orientation double covering map. 
\end{proposition}

\begin{proof}
Since a closed interval is Hausdorff, the condition $(2)$ implies the condition $(1)$. 
Denote by $p \colon \mathbb{S}^2 \to \mathbb{P}^2$ the orientation double cover. 
By Corollary~\ref{lifting_isotopic}, there is a lift $\widetilde{f}$ of $f$ along the orientation double covering map $p$. 
The lift $\widetilde{f}$ is a spherical homeomorphism isotopic to the identity. 
Since the lift $\widetilde{f}$ is isotopic to the identity, it is orientation-preserving. 
Proposition~\ref{lem43+} implies that $f$ is $R$-closed if and only if so is $\widetilde{f}$. 

Suppose that $f$ is $R$-closed. 
Then so is $\widetilde{f}$. 
Proposition~\ref{sphere} implies that the sphere consists of two fixed points and minimal sets each of which is a circloid. 
Since $\widetilde{f}$ is orientation-preserving, the projective plane $\mathbb{P}^2$ consists of one fixed point and minimal sets each of which is lifted into a circloid by the orientation double covering map. 

Suppose that the projective plane $\mathbb{P}^2$ consists of one fixed point and minimal sets each of which is lifted into a circloid by the orientation double covering map. 
Since $\widetilde{f}$ is orientation-preserving, the orientation double cover is a sphere which consists of two fixed points and minimal sets of $\widetilde{f}$ each of which is a circloid. 
Proposition~\ref{sphere} implies that $\widetilde{f}$ is $R$-closed. 
By  \cite[Theorem~6 and Theorem~7]{Ma}, the orbit class space of $\widetilde{f}$ is a closed interval. 
Let $H$ be the group generated by $\widetilde{f}$ and $G$ the group generated by $H$ and the covering  transformation group $\Gamma(p)$ of the orientation double covering map $p$. 
Then $G$ is isomorphic to $\Z \times \Z/2\Z$. 
The orbit class space $\mathbb{P}^2/\hat{f}$ of $f$ can be identified with the quotient space $(\mathbb{S}^2/\hat{H})/\Gamma(p) = \mathbb{S}^2/\hat{G}$ of the orbit class space $\mathbb{S}^2/\hat{H}$ of $\widetilde{f}$.  
Since the orbit class space of $\widetilde{f}$ is a closed interval, the quotient space $\mathbb{S}^2/\hat{H}$ by the $\Z/2\Z$-action of $\Gamma(p)$ is also a closed interval. 
This means that $\mathbb{P}^2/\hat{f}$ is a closed interval and so $f$ is $R$-closed. 
\end{proof}

We have the following characterization of $R$-closedness on a M\"obius band. 

\begin{corollary}
A homeomorphism $f$ isotopic to the identity on a M\"obius band is $R$-closed if and only if either it is periodic or each minimal set of $f$ is lifted into a finite disjoint union of essential circloids by the orientation double covering map.
\end{corollary}

\begin{proof}
Periodicity implies $R$-closedness. 
We may assume that $f$ is not periodic. 
By Corollary~\ref{lifting_isotopic}, there is a lift $\widetilde{f}$ on a closed annulus of $f$ along the orientation double covering map. 

Suppose that $f$ is $R$-closed. 
Proposition~\ref{lem43+} implies that the lift $\widetilde{f}$ by the orientation double covering map is $R$-closed. 
Since the boundary of a surface is invariant under any homeomorphism, the lift $\widetilde{f}$ is not minimal. 
By Lemma~\ref{lem:decomp_double}, taking the double of the annulus, the resulting surface is a torus and the double of the lift $\widetilde{f}$ is $R$-closed. 
Theorem~\ref{main:d} implies that each minimal set of the double of $\widetilde{f}$ and so of $\widetilde{f}$ is a finite disjoint union of essential circloids. 

Conversely, suppose that each minimal set of $f$ is lifted into a finite disjoint union of essential circloids by the orientation double covering map. 
The lift $\widetilde{f}$ is an annular homeomorphism isotopic to the identity. 
By taking the double of the annulus, the resulting surface is a torus and so each minimal set of the induced toral homeomorphism is a finite disjoint union of essential circloids. 
Then the orbit class space of $\widetilde{f}$ is either an interval or a circle, and so is Hausdorff. 
Theorem~\ref{main:d} implies the toral homeomorphism is $R$-closed and so is $\widetilde{f}$. 
Proposition~\ref{lem43+} implies that $f$ is $R$-closed. 
\end{proof}

The lifting condition in the previous lemmas is necessary. 
In fact, there is a homeomorphism isotopic to the identity on a M\"obius band $\mathbb{M}^2$ with a minimal set which is a circle but not a circloid. 
Indeed, let $\mathbb{A}^2 := (\R/2\Z) \times [0,1]$, and $\mathbb{M}^2 := ([0,1] \times [0,1])/(0,y)\sim(1,1-y)$ a M\"obius band. 
Consider an annular diffeomorphism $\widetilde{f}: \mathbb{A}^2 \to \mathbb{A}^2$ by $\widetilde{f}([x,y]) = ([x + \alpha, y])$ for some irrational number $\alpha$ and $p : \mathbb{A}^2 \to \mathbb{M}^2$ a covering map defined by $p([x,y]) = ([x,y])$ if $(x,y) \in [0,1]^2$ and $p([x,y]) = ([x-1, 1-y])$ if $(x,y) \in (1,2) \times [0,1]$.  
Then the induced homeomorphism $f \colon  \mathbb{M}^2 \to \mathbb{M}^2$ is $R$-closed but the circle $[[0,1] \times \{ 1/2 \}]$ is a minimal set but not annular because any neighborhood of the circle is non-orientable. 

Notice there is a minimal homeomorphism on a Klein bottle \cite{El,Pa}. 
We have the following characterization of $R$-closedness on a Klein bottle. 

\begin{corollary}
A homeomorphism isotopic to the identity on a Klein bottle is $R$-closed if and only if it satisfies one of the following statements:
\\
$(1)$ The homeomorphism $f$ is minimal.
\\
$(2)$ The homeomorphism $f$ is periodic.
\\
$(3)$ 
Each minimal set of $f$ is lifted into a finite disjoint union of essential circloids by the orientation double covering map. 
\end{corollary}

\begin{proof}
We may assume that $f$ is neither minimal nor periodic. 
By Corollary~\ref{lifting_isotopic}, there is a lift $\widetilde{f}$ of $f$ along the orientation double covering map which is a toral homeomorphism isotopic to the identity. 
Suppose that $f$ is $R$-closed. 
Proposition~\ref{lem43+} implies that $\widetilde{f}$ is $R$-closed. 
Theorem~\ref{main:d} implies that each minimal set of $\widetilde{f}$ is a finite disjoint union of essential circloids. 
Conversely, suppose that each minimal set of $f$ is lifted into a finite disjoint union of essential circloids by the orientation double covering map. 
Then each minimal set of $\widetilde{f}$ is a finite disjoint union of essential circloids. 
Theorem~\ref{main:d} implies that $\widetilde{f}$ is $R$-closed. 
Proposition~\ref{lem43+} implies that $f$ is $R$-closed. 
\end{proof}

The lifting condition in the previous lemmas is necessary. 
In fact, there is a homeomorphism isotopic to the identity on a Klein bottle with a minimal set which is a circle but not a circloid. 
Indeed, let $\mathbb{S}^1:=\R/\Z$ be a circle, $\mathbb{T}^2 := (\R/2\Z) \times \mathbb{S}^1$, and $\mathbb{K}^2 := ([0,1] \times \mathbb{S}^1)/(0,y)\sim(1,1-y)$ a Klein bottle. 
Consider a toral diffeomorphism $\widetilde{f}: \mathbb{T}^2 \to \mathbb{T}^2$ by $\widetilde{f}([x,y]) = ([x + \alpha, y])$ for some irrational number $\alpha$ and $p : \mathbb{T}^2 \to \mathbb{K}^2$ a covering map defined by $p([x,y]) = ([x,y])$ if $(x,y) \in [0,1]^2$ and $p([x,y]) = ([x-1, 1-y])$ if $(x,y) \in (1,2) \times [0,1]$.  
Then the induced homeomorphism $f \colon  \mathbb{K}^2 \to \mathbb{K}^2$ is $R$-closed but the circle $[[0,1] \times \{ 1/2 \}]$ is a minimal set but not annular because any neighborhood of the circle is non-orientable. 
We have the following characterization of $R$-closedness on a connected compact surface $S$ whose Euler number is negative. 

\begin{corollary}
Each homeomorphism isotopic to the identity on a connected compact surface $S$ whose Euler number is negative is periodic if and only if it is $R$-closed. 
\end{corollary}

\begin{proof}
By Corollary~\ref{lifting_isotopic}, we may assume that $S$ is non-orientable. 
Let $f$ be a homeomorphism isotopic to the identity on a connected compact surface $S$ whose Euler number is negative. 
As the orientable case, suppose that $f$ is periodic. 
Since the orbit class space of a periodic homeomorphism on a compact manifold is Hausdorff, the homeomorphism $f$ is $R$-closed. 
Conversely, suppose that $f$ is $R$-closed. 
By Corollary~\ref{lifting_isotopic}, there is a lift $\widetilde{f}$ of $f$ on the orientable double cover of $S$ which is an orientable compact surface. 
Then $\widetilde{f}$ is isotopic to the identity. 
Proposition~\ref{lem43+} implies that $\widetilde{f}$ is $R$-closed. 
Theorem~\ref{main:e} implies that $\widetilde{f}$ is periodic and so is $f$.  
\end{proof}

\section{Five possibilities of $R$-closed flows on compact $3$-manifolds}
In this section, we describe five possibilities of $R$-closed flows on connected compact $3$-manifolds. 
We obtain the following statements for flows on $3$-manifolds. 

\begin{lemma}\label{th:005}
Let $v$ be an $R$-closed flow on a connected compact $3$-manifold $M$. 
Then the union $U_C$ of suspensions of circloids is open and the quotient space $U_C/\hat{v}$ is a $1$-manifold. 
\end{lemma}

\begin{proof}
Let $\mathbb{A}$ be an open annulus transverse to $v$ and $C \subset \mathbb{A}$ a circloid whose saturation is a suspension $\mathcal{M}$ of a circloid with $C = \mathcal{M} \cap \mathbb{A}$ and $\mathcal{M} \cap \partial \mathbb{A} = \emptyset$. 
Then $\mathcal{M} \cap \mathrm{Sing}(v) = \emptyset$. 
Since $\mathcal{M}$ is $\hat{\F}_{v}$-invariant and connected, there are small connected neighborhoods $U_1, U_2 \subseteq \F_v(\mathbb{A})$ of $\mathcal{M}$ with $\overline{U_1}, \overline{U_2} \subseteq \F_v(\mathbb{A})$ such that the intersections $U_i \cap \mathbb{A}$ are deformation retractions of $\mathbb{A}$, and the first return map $f_v \colon  U_1 \cap \mathbb{A} \to U_2 \cap \mathbb{A}$ is well-defined. 
In particular, the intersections $U_i \cap \mathbb{A}$ are connected and annular. 
Note $\partial U_i \cap \mathcal{M} = \emptyset$. 
Since $v$ is $R$-closed and so $\hat{\F}_{v}$ is usc, there is an $\hat{\F}_{v}$-invariant connected open neighborhood $V \subseteq U_1 \cap U_2$ of $\mathcal{M}$. 

We claim that there is an $\hat{\F}_{v}$-invariant open connected neighborhood $W \subseteq V$ of $\mathcal{M}$ with $W = \F_v(W_\mathbb{A})$ such that the first return map $f_{W_\mathbb{A}} := f_v| \colon  W_\mathbb{A} \to W_\mathbb{A}$ is well-defined and homeomorphic, where $W_\mathbb{A} := W \cap \mathbb{A}$. 
Indeed, for a point $x \in C$, there is a small open connected neighborhood $W_x \subseteq V \cap \mathbb{A}$ of $x \in C$ in $\mathbb{A}$ such that the first return time map $t_x: W_x \to \R_{>0}$ with $f_v(y) = v(t_x(y), y) \in V \cap \mathbb{A}$ for any $y \in W_x$ is continuous.  
Since $C$ is compact, there are finitely many points $x_1, \ldots, x_n \in C$ such that $C \subset \bigcup_{i=1}^n W_{x_i}$. 
Then $W_C := \bigcup_{i=1}^n W_{x_i} \cap (\bigcup_{i=1}^n f_v(W_{x_i}))$ is an open neighborhood of $C$ in $\mathbb{A}$. 
Moreover the first return time map $t_C \colon  W_C \to \R_{>0}$ is well-defined and continuous.  
Since $v$ is $R$-closed, there is an $\F_{v}$-invariant open neighborhood $W \subseteq  \{ v([0, t_{C}(y)], y) \mid y \in W_C \}$ of $\mathcal{M}$. 
Put  $W_\mathbb{A} := W \cap \mathbb{A} \subseteq W_C$.  
Then $W_\mathbb{A} \subset U_1 \cap U_2 \cap \mathbb{A}$ is an open neighborhood of $C$ in $\mathbb{A}$ such that $\F_v(W_{\mathbb{A}}) = W$. 
Replacing $W$ with $\F_v(W_\mathbb{A})$ if necessary, we may assume that $W = \F_v(W_\mathbb{A})$. 
Since $W$ is $\F_{v}$-invariant, we have $f_v(W_\mathbb{A}) \subseteq W_\mathbb{A}$ and $f_v^{-1}(W_\mathbb{A}) \subseteq W_\mathbb{A}$. 
This means that $f_v(W_\mathbb{A}) = W_\mathbb{A}$ and so $f_v| \colon  W_\mathbb{A} \to W_\mathbb{A}$ is a homeomorphism and so that the claim holds. 

We claim that we may assume that $W_\mathbb{A}$ is connected. 
Indeed, take an open annular neighborhood $U \subseteq W_\mathbb{A}$ in $\mathbb{A}$ of the $f_v$-minimal set $C = \mathcal{M} \cap W_\mathbb{A}$. 
Since $U$ is arcwise-connected and $U \cap f_v^n(U) \neq \emptyset$ for any $n \in \Z$, the union $\bigcup_{n \in \Z} f_v^n(U) = \hat{\F}_{v}(U) \cap W_\mathbb{A}$ is arcwise-connected and so connected. 
Replacing $W$ with the open $\hat{\F}_{v}$-invariant connected neighborhood $\hat{\F}_{v}(U) \subseteq W$ of $\mathcal{M}$, we have that $\bigcup_{n \in \Z} f_v^n(U) = \hat{\F}_{v}(U) \cap \mathbb{A}$ is connected and so the claim holds. 

Let $B$ be the union of the connected components of $\mathbb{A} - W_\mathbb{A}$ which are contractible in $\mathbb{A}$. 
By $W_\mathbb{A} \subset U_1 \cap U_2 \cap \mathbb{A}$, since $U_1 \cap \mathbb{A}$ and $U_2 \cap \mathbb{A}$ are deformation retractions of $\mathbb{A}$, we have $B \subset U_1 \cap U_2 \cap \mathbb{A}$. 
Then it suffices to show that $B \cup W_\mathbb{A}$ is an open annulus and $\F_v(B \cup W_\mathbb{A})$ is an $\hat{\F}_{v}$-invariant connected open neighborhood of $\mathcal{M}$ which consists of suspensions of circloids such that the orbit class space $\F_v(B \sqcup W_\mathbb{A})/\hat{v}$ is an open interval. 
%
Recall the filling $\mathrm{Fill}_\mathbb{A}(W_\mathbb{A})$ as follows: $p \in \mathrm{Fill}_\mathbb{A}(W_\mathbb{A})$ if either $p \in W_\mathbb{A}$ or there is an open disk in $\mathbb{A}$ which  contains $p \in \mathbb{A}$ and whose boundary is contained in $W_\mathbb{A}$. 
Since $W_{\mathbb{A}}$ is open and any connected components of $\partial W_{\mathbb{A}}$ which are not contractible in $\mathbb{A}$ are parallel to $\mathbb{A}$, we have that $p \in \mathrm{Fill}_\mathbb{A}(W_\mathbb{A})$ if and only if there is a simple closed curve in $W_\mathbb{A}$ which bounds  a disk in $\mathbb{A}$ containing $p \in \mathbb{A}$. 

We claim that $\mathrm{Fill}_\mathbb{A}(W_\mathbb{A})$ is an open annulus and that $f' := f_v| \colon  \mathrm{Fill}_\mathbb{A}(W_\mathbb{A}) \to \mathrm{Fill}_\mathbb{A}(W_\mathbb{A})$ is a homeomorphism. 
Indeed, 
the facts that $U_i \cap \mathbb{A}$ are annular imply that $U_i = \mathrm{Fill}_\mathbb{A}(U_i)$ and so $C \subset \mathrm{Fill}_\mathbb{A}(W_\mathbb{A}) \subseteq U_1 \cap U_2 = \mathrm{dom}(f_v) \cap \mathrm{im}(f_v)$. 
Since each point of $B$ is bounded by a simple closed curve in $W_\mathbb{A}$, this implies that the filling $\mathrm{Fill}_\mathbb{A}(W_\mathbb{A}) = B \sqcup W_\mathbb{A}$ is annular such that $C \subset f_v(\mathrm{Fill}_\mathbb{A}(W_\mathbb{A})) = \mathrm{Fill}_\mathbb{A}(W_\mathbb{A}) \cup f_v(B) \subseteq f_v(U_1) \subseteq U_2 \subseteq \mathbb{A}$.
By $\partial B \cup f_v(\partial B) \cup f_v^{-1}(\partial B) \subseteq \partial W_\mathbb{A}$ and $f_v(W_\mathbb{A}) = f_v^{-1}(W_\mathbb{A}) = W_\mathbb{A}$, since $B$ is the union of the connected components of $\mathbb{A} - W_\mathbb{A}$ which are contractible in $\mathbb{A}$, the image $f(B)$ and the inverse image $f_v^{-1}(B)$ are subsets of the union of the connected components of $U_1 \cap U_2 \cap (\mathbb{A} - W_\mathbb{A})$ which are contractible in $\mathbb{A}$. 
This means that $f_v(B) \cup f_v^{-1}(B) \subseteq B$. 
The injectivity of $f_v$ implies that $f_v(B) = B$, and so that $f' := f_v| \colon  \mathrm{Fill}_\mathbb{A}(W_\mathbb{A}) \to \mathrm{Fill}_\mathbb{A}(W_\mathbb{A})$ is a homeomorphism. 
Write $W_F := \F_v(\mathrm{Fill}_\mathbb{A}(W_\mathbb{A})) = \hat{\F}_v(\mathrm{Fill}_\mathbb{A}(W_\mathbb{A}))$. 
From $\F_v(W_{\mathbb{A}}) = W$, we have $\mathbb{A} \cap W_F = \mathrm{Fill}_\mathbb{A}(W_\mathbb{A})$. 
Since $C$ is a circloid and $\mathbb{A}$ is the annular neighborhood, we obtain that $\mathrm{Fill}_\mathbb{A}(W_\mathbb{A})$ is an open annulus. 
The $R$-closedness of $v$ implies that $M/\hat{v}$ is Hausdorff and so is $\mathrm{Fill}_\mathbb{A}(W_\mathbb{A})/\hat{\F}_{f'} = \mathrm{Fill}_\mathbb{A}(W_\mathbb{A})/\hat{\F}_{f_v}$. 
Since the boundary of a saturation is invariant, by$\mathbb{A} \cap \partial W_F = \partial_{\mathbb{A}} (\mathrm{Fill}_\mathbb{A}(W_\mathbb{A}))$, we have $\hat{\F}_{v}(\partial_{\mathbb{A}} (\mathrm{Fill}_\mathbb{A}(W_\mathbb{A}))) \subseteq \partial (\hat{\F}_{v}(\mathrm{Fill}_\mathbb{A}(W_\mathbb{A}))) = \partial W_F$, where $\partial_{\mathbb{A}} F$ is the boundary of a subset $F \subseteq \mathbb{A}$ in $\mathbb{A}$.  

We claim that $\hat{\F}_{v}(\partial_{\mathbb{A}} (\mathrm{Fill}_\mathbb{A}(W_\mathbb{A}))) = \partial W_F$.  
Indeed, otherwise there is a point $x \in \partial W_F - \hat{\F}_{v}(\partial_{\mathbb{A}} (\mathrm{Fill}_\mathbb{A}(W_\mathbb{A})))$. 
Since the saturation $\hat{\F}_{v}(\partial_{\mathbb{A}} (\mathrm{Fill}_\mathbb{A}(W_\mathbb{A})))$ is compact, there are finitely many points $x_1, \ldots, x_k \in  \hat{\F}_{v}(\partial_{\mathbb{A}} (\mathrm{Fill}_\mathbb{A}(W_\mathbb{A})))$ and their $\hat{\F}_{v}$-invariant neighborhoods $U_i$ of $x_i$ with $x \notin \overline{U_i}$ such that $\hat{\F}_{v}(\partial_{\mathbb{A}} (\mathrm{Fill}_\mathbb{A}(W_\mathbb{A}))) \subset \bigcup_{i=1}^k U_i$. 
Then $W_F \cup \bigcup_{i=1}^k \overline{U_i}$ is a closed $\hat{\F}_{v}$-invariant neighborhood of $W_F$ which does not intersect $x$. 
This implies that $x \notin \overline{W_F}$, which contradicts to the choice of $x$.

By the two-point compactification of $\mathrm{Fill}_\mathbb{A}(W_\mathbb{A})$, we obtain the resulting sphere $S$ and the resulting homeomorphism $f_{S}$ from $f'$ with the two periodic points which are the new points. 
By the construction, 
the class space $S/\hat{\F}_{f_S}$ is the two-point compactification of $\mathrm{Fill}_\mathbb{A}(W_\mathbb{A})/\hat{\F}_{f'}$. 
Put $f'_S := f_S$ if $f_S$ is orientation-preserving and $f'_S := f_S^2$ if $f_S$ is orientation-reversing. 
Since $M$ is normal, the $R$-closedness of $v$ implies that $S/\hat{\F}_{f'_S}$ is Hausdorff and so $f'_S$ is $R$-closed. 
By \cite[Corollary~2.6]{Y2}, we obtain that $S$ consists of two fixed points and circloids, and so that the orbit class space $S/\hat{\F}_{f'_S}$ is a closed interval. 
Because the quotient space of a $1$-manifold by a finite group-action is a $1$-orbifold, which is homeomorphic to a $1$-manifold with or without boundary. 
Since $f'_S$ is either $f_S$ or $f_S^2$, the quotient space $S/\hat{\F}_{f_S}$ is a closed interval and so is $\overline{\mathrm{Fill}_\mathbb{A}(W_\mathbb{A})}/\hat{\F}_{f_v}$. 
Therefore the orbit class space $\mathrm{Fill}_\mathbb{A}(W_\mathbb{A})/\hat{\F}_{f_S} = \mathrm{Fill}_\mathbb{A}(W_\mathbb{A})/\hat{\F}_{f_v} = \F_v(\mathrm{Fill}_\mathbb{A}(W_\mathbb{A}))/\hat{v} = \F_v(B \sqcup W_\mathbb{A})/\hat{v}$ is an open interval and the open neighborhood $\F_v(B \sqcup W_\mathbb{A}) = \F_v(\mathrm{Fill}_\mathbb{A}(W_\mathbb{A})) = W_F$ of $\mathcal{M}$ consists of suspensions of circloids. 

Thus the quotient space $U_C/\hat{v} \subseteq M/\hat{v}$ of the union of suspensions of circloids is an union of open intervals and is locally Euclidean. 
Since $M/\hat{v}$ is Hausdorff, the quotient space $U_C/\hat{v} \subseteq M/\hat{v}$ is Hausdorff and so a $1$-manifold.  
\end{proof}


We have the following tameness of transversality of minimal sets. 

\begin{lemma}\label{lem52} 
Let $v$ be an $R$-closed flow on a connected compact $3$-manifold $M$. 
Suppose that each two-dimensional minimal set is a suspension of a circloid. 
If there is a suspension of a circloid, then the orbit class space $M/\hat{v}$ of $M$ is a closed interval or a circle such that there are at most two orbits each of which is not the suspension of a circloid. 
\end{lemma}

\begin{proof} 
By Lemma~\ref{th:005}, the union $U$ of suspensions of circloids is open and the quotient space of it is a $1$-manifold. 
The Hausdorff property of $M/\hat{v}$ implies that $\overline{U}/\hat{v}$ is either a closed interval or a circle. 
The $R$-closedness implies that $\partial U$ consists of at most two orbit closures each of whose dimensions of $\partial U$ is at most one.   
By Lemma~\ref{lem062}, we have that $\overline{U} = M$ and so that $M/\hat{v} $ is a closed interval or a circle. 
\end{proof}

A subset $C$ in a topological space $X$ is separating if the complement $X - C$ is disconnected. 
The following statement is obtained by a similar argument as the proof of Lemma~\ref{lem55}. 

\begin{lemma}\label{lem54}
Let $v$ be an $R$-closed flow on a connected compact $3$-manifold $M$. 
Each one-dimensional minimal set is a periodic orbit. 
\end{lemma}

\begin{proof} 
Assume that there is a one-dimensional minimal set $\mathcal{M}$ which is not a periodic orbit. 
Since $\mathcal{M} \cap \mathrm{Sing}(v) = \emptyset$, there is a neighborhood of $\mathcal{M}$ without singular points. 
By the $R$-closedness of $v$, there is an open $\hat{\F}_{v}$-invariant neighborhood $U$ of $\mathcal{M}$ without singular points. 
The minimality of $\mathcal{M}$ implies that, for any point $x \in \mathcal{M}$ and any transverse closed disk $T_0 \subset U$ to $v$ with $x \in T_0$, the point $x$ is not isolated in the intersection $\mathcal{C}_0 := T_0 \cap \mathcal{M}$. 
Then the intersection $\mathcal{C}_0 = T_0 \cap \mathcal{M}$ is a zero-dimensional compact metrizable space. 
Since a compact Hausdorff space is zero-dimensional if and only if it is totally disconnected (cf.  \cite[Proposition~3.1.7, p.136]{AT2008}), the intersection $\mathcal{C}_0$ is totally disconnected. 

We claim that there is a transverse open disk $T \subset T_0 \subset U$ to $v$ such that $\mathcal{C} := T \cap \mathcal{M}$ is a Cantor set with $\mathcal{C} \cap \partial T = \emptyset$ and $\F_v(\mathcal{C}) = \hat{\F}_v(\mathcal{C}) = \mathcal{M}$. 
Indeed, since a zero-dimensional space has a basis consisting of open and closed subsets, for any point $x \in \mathcal{C}$, there is an open disk $T$ containing $x$ with $\overline{T} \subset \mathrm{int} T_0$ and $\mathcal{C} \cap \partial T = \emptyset$. 
Then the intersection $\mathcal{C} := T \cap \mathcal{M} = \overline{T}  \cap \mathcal{M}$ is a nonempty perfect zero-dimensional compact metrizable subset. 
Since a Cantor set is a compact metrizable perfect totally disconnected space, the intersection $\mathcal{C} = T \cap \mathcal{M}$ is a Cantor set. 
The minimality of $\mathcal{M}$ implies that $\overline{O(y)} = \mathcal{M}$ for any point $y \in \mathcal{M}$, and so that $\emptyset \neq O(y) \cap T =  O(y) \cap \mathcal{C}$ for any point $y \in \mathcal{M}$. 
Therefore $\hat{\F}_v(\mathcal{C}) = \F_v(\mathcal{C}) = \mathcal{M}$. 
This completes the proof of the claim. 

%
%
For any $x \in \mathcal{C}$, since a Cantor set has a basis consisting of open and closed subsets, there is an open disk $V_x \subseteq T$ with $\overline{V_x} \subset T$ and $\mathcal{C} \cap \partial V_x = \emptyset$ such that  the first return time map $t_x \colon \overline{V_x} \to \R_{>0}$ to $T$ with $v(t_x(y), y) \in T$ for any $y \in \overline{V_x}$ is continuous. 
Then the saturation $\F_v(V_x)$ of the open transverse disk $V_x$ is open and $\F_v(V_x) = \hat{\F}_v(V_x)$ because of Lemma~\ref{lem02}. 
Since $\mathcal{C}$ is compact, there are finitely many points $x_1, \ldots , x_k \in \mathcal{C}$ such that $\mathcal{C} \subset \bigcup_{i =1}^k V_{x_i} \subset \bigcup_{i =1}^k \overline{V_{x_i}} \subset T$. 
Put $D := \hat{\F}_v(\bigcup_{i =1}^k \partial V_{x_i})$. 
The $R$-closedness of $v$ implies that $D$ is closed. 

We claim that $\mathcal{M}$ and $D$ are disjoint invariant closed subsets. 
Indeed, recall that the filling $\mathrm{Fill}_{T}(B)$ of a subset $B \subseteq T$ in the open disk $T$ as follows: $y \in \mathrm{Fill}_{T}(B)$ if either $y \in B$ or there is an open disk in $T$ containing $y$ whose boundary is contained in $B$. 
Denote by $F_\mathcal{C} := \mathrm{Fill}_{T}(\bigcup_{i =1}^k \overline{V_{x_i}}) \subset T$ the filling of the union $\bigcup_{i =1}^k \overline{V_{x_i}}$. 
Then the first return map $f_v: \bigcup_{i=1}^k \overline{V_{x_i}} \to T$ is continuous. 
Since $\mathcal{C} = T \cap \mathcal{M} = T \cap \mathcal{M}$, $\mathcal{M} \cap \partial T = \emptyset$, and $\F_v(\mathcal{C}) = \mathcal{M}$, the minimality of $\mathcal{M}$ implies that $f_v(\mathcal{C}) = \mathcal{C}$. 
Since $V_{x_i}$ are open disks in $T$, the filling $F_\mathcal{C}$ consists of finitely many connected components. 
By definition, we have $\mathcal{M} \cap \bigcup_{i=1}^k \partial V_{x_i} = \mathcal{C} \cap \bigcup_{i=1}^k \partial V_{x_i} = \emptyset$. 
The invariance of $\mathcal{M}$ implies $\mathcal{M} \cap D = \emptyset$. 
This completes the claim. 

By $f_v(\mathcal{C}) = \mathcal{C} \subset \bigcup_{i =1}^k V_{x_i} \subset \mathrm{dom}(f_v) \subseteq T$, we obtain $f_v(\mathcal{C}) = \mathcal{C} \subset (\bigcup_{i =1}^k V_{x_i} ) \cap (\bigcup_{i =1}^k f_v(V_{x_i}) )$ and so $f_v(\mathcal{C}) = \mathcal{C} \subset \mathrm{int}\left(\bigcup_{i =1}^k \{ v([0, t_{x_i}(y)], y) \mid y \in V_{x_i} \}\right)$. 
Therefore there are disjoint open connected $\hat{\F}_{v}$-invariant neighborhoods $U_{D}$ and $U_{\mathcal{M}} \subseteq \bigcup_{i =1}^k \{ v([0, t_{x_i}(y)], y) \mid y \in V_{x_i} \} \subseteq \bigcup_{i =1}^k \F_v(V_{x_i}) = \bigcup_{i =1}^k \hat{\F}_v(V_{x_i})$ of $D$ and $\mathcal{M}$ respectively. 
%
Put $V_{\mathcal{M}} := T \cap U_{\mathcal{M}} \subseteq \mathrm{dom}(f_v) \cap \mathrm{im}(f_v)$. 
Then $V_{\mathcal{M}}$ is open in $T$. 
Since $\mathcal{C} = T \cap \mathcal{M}$, we have $\mathcal{C}  = T \cap \mathcal{M} \subset T \cap U_{\mathcal{M}} = V_{\mathcal{M}}$ and so $V_{\mathcal{M}} \subseteq \bigcup_{i =1}^k V_{x_i}$. 
Then $V_{\mathcal{M}} =   \bigcup_{i =1}^k V_{\mathcal{M}} \cap V_{x_i}$. 

We claim that we may assume that $\F_v(V_{\mathcal{M}}) = U_{\mathcal{M}}$. 
Indeed, taking an open connected \nbd $V'_{\mathcal{M}} \subseteq V_{\mathcal{M}}$ of $\mathcal{C}$ in $T$ with $\overline{V'_{\mathcal{M}}} \subset V_{\mathcal{M}}$ and replacing $U_{\mathcal{M}}$ by $\hat{\F}_v(V'_{\mathcal{M}})$ if necessary, we have $\hat{\F}_v(V_{\mathcal{M}}) = \F_v(V_{\mathcal{M}}) = U_{\mathcal{M}}$. 

%
Since $V_{\mathcal{M}} \cap \bigcup_{i =1}^k \partial V_{x_i} \subseteq U_{\mathcal{M}} \cap U_{D} = \emptyset$, any connected component of $V_{\mathcal{M}}$ is contained in $V_{x_j} \setminus \bigcup_{i =1}^k \partial V_{x_i}$ for some $j$. 
This means that the intersection $V_{\mathcal{M}} \cap V_{x_i}$ is the union of connected components of $V_{\mathcal{M}}$ contained in $V_{x_i}$. 
By the openness of $V_{\mathcal{M}} \cap V_{x_i}$ in the open disk $V_{x_i}$, each connected component of the intersection $V_{\mathcal{M}} \cap V_{x_i}$ is a punctured open disk in $V_{x_i} \setminus \bigcup_{j =1}^k \partial V_{x_j}$.
Then $\partial_T V_{\mathcal{M}} = \bigcup_{i =1}^k \partial_T (V_{\mathcal{M}} \cap V_{x_i})$. 
Moreover, $ p \in \mathrm{Fill}_{V_{x_i}}(V_{\mathcal{M}} \cap V_{x_i})$ if and only if there is a simple closed curve in $V_{\mathcal{M}} \cap V_{x_i}$ which bounds a disk in $V_{x_i}$ containing $p$. 
Put $V_i := \mathrm{Fill}_{V_{x_i}}(V_{\mathcal{M}} \cap V_{x_i}) \subseteq V_{x_i}$ and $V_{\mathcal{C}} := \bigcup_{i =1}^k V_i \subseteq \bigcup_{i =1}^k V_{x_i} \subseteq \bigcup_{i =1}^k \overline{V_{x_i}} \subseteq T$. 
Then the subsets $V_i$ and $V_{\mathcal{C}}$ are open in $T$ such that $\overline{V_{\mathcal{C}}} \subset T$. 
Moreover, the filling $V_i = \mathrm{Fill}_{V_{x_i}}(V_{\mathcal{M}} \cap V_{x_i})$ is the filling in the transverse open disk $V_{x_i}$ of the union of connected components of $V_{\mathcal{M}}$ contained in $V_{x_i}$ and so  
 is a disjoint union of open disks in $V_{x_i} \subset T$ such that $\partial_T V_i \subseteq \partial_T (V_{\mathcal{M}} \cap V_{x_i}) \subseteq \partial_T V_{\mathcal{M}}$. 
Therefore the union $V_{\mathcal{C}}$ is a disjoint union of open disks such that $\partial_T V_{\mathcal{C}} \subseteq \bigcup_{i=1}^k \partial_T (V_{\mathcal{M}} \cap V_{x_i}) = \partial_T V_{\mathcal{M}}$. 
%
%
Put $U_{\mathcal{C}} := \F_v(V_{\mathcal{C}})$. 
Proposition~\ref{corollary:cl_sat} implies that $\hat{\F}_v(\overline{V_{\mathcal{C}}}) = \overline{\F_v(V_{\mathcal{C}})} = \overline{U_{\mathcal{C}}}$. 
The transversality of $T$ implies that the saturation $U_{\mathcal{C}} = \hat{\F}_v(V_{\mathcal{C}})$ is open and so that $\partial U_{\mathcal{C}} = \overline{U_{\mathcal{C}}} - U_{\mathcal{C}} = \hat{\F}_v(\overline{V_{\mathcal{C}}}) - \hat{\F}_v(V_{\mathcal{C}}) = \hat{\F}_v(\partial V_{\mathcal{C}})$.
By the $R$-closedness of $v$, since the boundary $\partial U_{\mathcal{C}}$ is invariant and closed, the boundary $\partial U_{\mathcal{C}}$ consists of minimal sets of $v$. 
Since any connected component of $V_{\mathcal{M}}$ is contained in $V_{x_j} \setminus \bigcup_{i =1}^k \partial V_{x_i}$ for some $j$, so is its filling. 
In other words, any connected component of $V_{\mathcal{C}}$ is an open disk in $V_{x_j} \setminus \bigcup_{i =1}^k \partial V_{x_i}$ for some $j$. 
%

We claim that $V_{\mathcal{C}} = \bigcup_{i =1}^k \mathrm{Fill}_{V_{x_i}}(V_{\mathcal{M}} \cap V_{x_i})$ is $f_v$-invariant. 
Indeed, we have $U_{\mathcal{M}} = \F_v(V_{\mathcal{C}}) \subseteq  \F_v(\bigcup_{i =1}^k V_{x_i})$. 
Since $U_{\mathcal{M}}$ is $\hat{\F}_{v}$-invariant, the intersection $V_{\mathcal{M}} = T \cap U_{\mathcal{M}} \subseteq \mathrm{dom}(f_v) \cap \mathrm{im}(f_v) \subseteq \bigcup_{i =1}^k \overline{V_{x_i}}$ is $f_v$-invariant. 
Fix a simple closed curve $\gamma$ in $V_{\mathcal{M}} \subseteq V_{\mathcal{C}}$. 
Then $\gamma \subseteq V_{\mathcal{M}} \cap V_{x_j}$ for some $j$. 
Since $f_v(V_{\mathcal{M}}) = V_{\mathcal{M}}$, the image $f_v(\gamma)$ is contained in $V_{\mathcal{M}}$. 
By $V_{\mathcal{M}} \cap (\bigcup_{i =1}^k \partial V_{x_i}) = \emptyset$, we have $f_v(\gamma) \subset V_{\mathcal{M}} \cap V_{x_k}$ for some $k$.
Since $V_{x_i}$ is an open disk, the fillings $\mathrm{Fill}_{V_{x_j}}(\gamma) \subset V_{x_j} \subseteq \mathrm{dom}(f_v)$ and $\mathrm{Fill}_{V_{x_k}}(f_v(\gamma)) \subset V_{x_k}$ are closed disks in $T$. 
Then the image $f_v(\mathrm{Fill}_{V_{x_j}}(\gamma)) \subset \mathrm{im}(f_v) \subseteq T$ is a closed disk in the open disk $T$. 
Since two closed disks $\mathrm{Fill}_{V_{x_k}}(f_v(\gamma))$ and $f_v(\mathrm{Fill}_{V_{x_j}}(\gamma))$ contained in the open disk $T$ such that $\partial (\mathrm{Fill}_{V_{x_k}}(f_v(\gamma))) = f_v(\gamma) = \partial (f_v(\mathrm{Fill}_{V_{x_j}}(\gamma)))$, the Jordan curve theorem implies that  $\mathrm{Fill}_{V_{x_k}}(f_v(\gamma)) = f_v(\mathrm{Fill}_{V_{x_j}}(\gamma))$ and so $\mathrm{Fill}_{V_{x_k}}(f_v(V_{\mathcal{M}} \cap V_{x_k}))) = f_v(\mathrm{Fill}_{V_{x_j}}(V_{\mathcal{M}} \cap V_{x_j})))$. 
Since $V_{\mathcal{M}} = \bigcup_{i =1}^k V_{\mathcal{M}} \cap V_{x_i}$ is $f_v$-invariant, this means that $f_v(V_{\mathcal{C}}) = f_v(\bigcup_{i =1}^k \mathrm{Fill}_{V_{x_i}}(V_{\mathcal{M}} \cap V_{x_i})) = \bigcup_{i =1}^k f_v(\mathrm{Fill}_{V_{x_i}}(V_{\mathcal{M}} \cap V_{x_i})) = \bigcup_{i =1}^k \mathrm{Fill}_{V_{x_i}}(V_{\mathcal{M}} \cap V_{x_i}) = V_{\mathcal{C}}$. 

We claim that $f_v(\overline{V_{\mathcal{C}}}) = \overline{V_{\mathcal{C}}} \subset T$. 
Indeed, since the domain $\bigcup_{i=1}^k \overline{V_{x_i}}$ is compact and the codomain $T$ is Hausdorff, the image $f_v(\overline{V_{\mathcal{C}}})$ of a compact subset is closed. 
By $V_{\mathcal{C}} =  f_v(V_{\mathcal{C}}) \subset f_v(\overline{V_{\mathcal{C}}})$, we have $\overline{V_{\mathcal{C}}} =  \overline{f_v(V_{\mathcal{C}})} \subseteq f_v(\overline{V_{\mathcal{C}}})$. 
Conversely, the injectivity of $f_v$ implies that $V_{\mathcal{C}} = f_v^{-1}(V_{\mathcal{C}})$. 
The continuity of $f_v$ implies $\overline{V_{\mathcal{C}}} \subseteq f_v^{-1}(\overline{V_{\mathcal{C}}})$. 
By the injectivity of $f_v$, we have $f_v(\overline{V_{\mathcal{C}}}) \subseteq f_v(f_v^{-1}(\overline{V_{\mathcal{C}}})) = \overline{V_{\mathcal{C}}}$. 

Moreover, by the invariance of $\overline{V_{\mathcal{C}}}$ for the first return map $f_v$, the restriction $f_{V_\mathcal{M}} := f_v| \colon  \overline{V_{\mathcal{C}}} \to \overline{V_{\mathcal{C}}}$ is a homeomorphism because the restriction and the inverse map are well-defined and continuous injections.  
Thus the boundary $\partial_T V_{\mathcal{C}}$ is $f_v$-invariant such that $T \cap \partial (U_{\mathcal{C}} ) = T \cap \partial (\F_v(V_{\mathcal{C}})) = T \cap \F_v(\partial_T V_{\mathcal{C}}) = \partial_T V_{\mathcal{C}}$. 
Let $D_{x_1} \subset V_{x_1} \cap V_{\mathcal{C}}$ be the connected component of $V_{\mathcal{C}}$ containing $x_1 \in \mathcal{C}$. 
As mention above, the connected component $D_{x_1}$ is an open disk in $V_{x_1}$. 
The minimality of $\mathcal{M}$ implies that $\F_v(D_{x_1}) = \F_v(V_{\mathcal{C}}) = U_{\mathcal{C}}$  
 and $\hat{\F}_v(\partial_T D_{x_1}) = \hat{\F}_v(\partial_T V_{\mathcal{C}}) = \hat{\F}_v(T \cap \partial U_{\mathcal{C}}) = \F_v(T \cap \partial U_{\mathcal{C}}) \subseteq \partial U_\mathcal{C} = \partial (\F_v(V_\mathcal{C}))$. 
 
We claim that  $\hat{\F}_v(\partial_T V_{\mathcal{C}}) = \partial U_\mathcal{C}$. 
Indeed, otherwise there is a point $y \in \partial U_\mathcal{C} - \hat{\F}_v(\partial_T V_{\mathcal{C}})$. 
Since $\hat{\F}_v(\partial_T V_{\mathcal{C}})$ is compact, there are finitely many points $y_1, \ldots, y_l \in \hat{\F}_v(\partial_T V_{\mathcal{C}})$ and their $\hat{v}$-invariant neighborhoods $U_i$ of $y_i$ such that $\hat{\F}_v(\partial_T V_{\mathcal{C}}) \subset \bigcup_{i=1}^l U_i$ and $y \notin \overline{U_i}$. 
The union $\hat{\F}_v(V_{\mathcal{C}}) \cup \bigcup_{i=1}^l \overline{U_i}$ is a closed $\hat{v}$-invariant neighborhood of $U_{\mathcal{C}}$ which does not intersect $y$. 
This implies that $y \notin \overline{U_{\mathcal{C}}}$, which contradicts to the choice of $y$. 

Threrefore $\hat{\F}_v(\overline{D_{x_1}}) = \overline{U_\mathcal{C}}$. 
By $D_{x_1} \subset V_{x_1} \cap V_{\mathcal{C}}$, the first return map $f_D \colon \overline{D_{x_1}} \to \overline{D_{x_1}}$ on $\overline{D_{x_1}}$  induced by the homeomorphism $f_{V_\mathcal{M}} \colon \overline{V_{\mathcal{C}}} \to \overline{V_{\mathcal{C}}}$ and the inverse map are well-defined and continuous injections. 
This implies that $f_D$ is a homeomorphism and so is the restriction $f_D|_{D_{x_1}}$. 
By Theorem~\ref{th:usc}, the $R$-closedness of $v$ implies that the class space  $\overline{D_{x_1}}/\hat{f}_D$ which is homeomorphic to $\F_v(\overline{D_{x_1}})/\hat{v}$ is Hausdorff.
Therefore $f_D$ is $R$-closed. 
By the one-point compactification of $D_{x_1}$, the resulting space $S := (D_{x_1}/\hat{v}) \sqcup \{ \infty \}$ is a sphere. 
Moreover, the sphere $S$ can be obtained by collapsing the boundary $\partial_T D_{x_1}$ into a singleton. 
The induced homeomorphism $f_S \colon S \to S$ can be obtain from the restriction $f_D|_{D_{x_1}} \colon D_{x_1} \to D_{x_1}$ by adding a new fixed point. 
In addition, the homeomorphism $f_S$ also can be obtained from $f_D$ by collapsing the boundary $\partial_T D_{x_1}$ into a new fixed point. 
%
Since the boundary $\partial U_{\mathcal{C}}$ is closed $\F_v$-invariant, the orbit class space $S/\hat{f}_S$ also can be identified with the resulting surface of the orbit class space $\overline{D_{x_1}}/\hat{f}_D$ collapsing $\partial_T D_{x_1}/\hat{f}_D$ into a point. 
Then the orbit class space $S/\hat{f}_S$ is homeomorphic to the resulting space of $\overline{U_{\mathcal{C}}}/\hat{v}$ collapsing the quotient $\partial U_{\mathcal{C}}/\hat{v}$ of the boundary $\partial U_{\mathcal{C}}$ into a point. 
Thus the normality of $M$ implies that the orbit class space $S/\hat{f}_S$ is Hausdorff. 
Theorem~\ref{th:usc} implies that $f_S$ is an $R$-closed homeomorphism on the sphere $S$ with a Cantor minimal set. 
This contradicts to  \cite[Corollary~2.6]{Y2}. 
Thus each one-dimensional minimal set is a periodic orbit. 
\end{proof}

Under a non-existence of wild two-dimensional minimal sets, the existence of three periodic orbits implies that open denseness of periodic orbits as follows.  

\begin{lemma}\label{lem:006} 
Let $v$ be an $R$-closed flow on a connected compact $3$-manifold $M$. Suppose that each two-dimensional minimal set is a suspension of a circloid. 
If there are at least three distinct periodic orbits, then the periodic point set $\mathrm{Per}(v)$ is open dense, the quotient $\mathrm{Per}(v)/v$ is an open connected surface, the end set of $\mathrm{Per}(v)/v$ is zero-dimensional and homeomorphic to $\mathrm{Sing}(v)$, and $M = \mathrm{Sing}(v) \sqcup \mathrm{Per}(v)$. 
\end{lemma}

\begin{proof} 
If there is a three-dimensional minimal set, then the $R$-closedness implies that it has no boundaries and so equals  $M$. 
Thus there are no three-dimensional minimal sets. 
From Lemma~\ref{lem52}, there are no two-dimensional orbit closures.  
Lemma~\ref{lem54} implies that zero (resp. one) dimensional minimal sets are singular points (resp. periodic orbits), and so that $M = \mathrm{Sing}(v) \sqcup \mathrm{Per}(v)$. 
Since $\mathrm{Sing}(v)$ is closed, we have that $\mathrm{Per}(v)$ is open. 
By \cite[Theorem~3.12]{D3}, the quotient space $\mathrm{Per}(v)/\hat{v} = \mathrm{Per}(v)/v$ of $\mathrm{Per}(v)$ is an open connected surface. 
Then each end of $\mathrm{Per}(v)/\hat{v}$ corresponds to a connected component of $\partial \mathrm{Sing}(v)$. 
Since $M/\hat{v} $ is compact Hausdorff and so normal, the dimension of $\partial(\mathrm{Per}(v)/\hat{v})$ is at most one. 
Since $\partial(\mathrm{Per}(v)/\hat{v})$ consists of singular points, the Hausdorff separation axiom for $M/\hat{v}$ implies that any connected components of $\partial \mathrm{Per}(v) = \partial \mathrm{Sing}(v)$ is a singular point. 
By \cite[Theorem 3]{R}, ant open surface without boundary is homeomorphic to the resulting surface from a sphere by removing a closed totally disconnected subset. 
This means that $\mathrm{Sing}(v) = \partial \mathrm{Sing}(v)$ corresponds to the end set of $\mathrm{Per}(v)/\hat{v}$ which is homeomorphic to a closed totally disconnected subset in a sphere. 
Since a compact Hausdorff space is zero-dimensional if and only if it is totally disconnected (cf.  \cite[Proposition~3.1.7, p.136]{AT2008}), the singular point set $\mathrm{Sing}(v)$ is zero-dimensional. 
By Lemma~\ref{lem062}, we obtain $\overline{\mathrm{Per}(v)} = M$. 
\end{proof}

\subsection{Proof of Theorem~\ref{main:f}}
%
%
Let $v$ be a flow on a connected compact $3$-manifold $M$.  
From Theorem~\ref{th:usc}, one of cases $(1)$--$(3)$ implies that 
the class space is Hausdorff and so that $v$ is $R$-closed. 
Suppose that the case $(4)$ holds. 
Then the quotient $\mathrm{Per}(v)/v$ is an open connected surface whose end set is zero-dimensional and homeomorphic to $\mathrm{Sing}(v)$. 
Since $M = \mathrm{Sing}(v) \sqcup \mathrm{Per}(v)$, this implies that the class space $M/\hat{v} = M/v$ is Hausdorff and so $v$ is $R$-closed. 

Conversely, suppose that $v$ is $R$-closed. 
We may assume that $v$ is non-trivial (i.e. neither minimal nor identical). 
If there is a three-dimensional minimal set, then the $R$-closedness implies that it has no boundaries and so equals  $M$. 
Thus there are no three-dimensional minimal sets. 
Suppose that each two-dimensional minimal set is a suspension of a circloid. 
Since $v$ is non-trivial, there is a minimal set $\mathcal{M}$ which is not a singular point. 
If the minimal set $\mathcal{M}$ is two-dimensional, then Lemma~\ref{lem52} implies (3). 
Thus we may assume that there are no two-dimensional minimal sets. 
By Lemma~\ref{lem54}, each one-dimensional minimal set is a periodic orbit and so $\mathcal{M}$ is a periodic orbit. 
Because the singular point set is closed, this has a neighborhood without singular points. 
Since there are neither three-dimensional nor two-dimensional minimal sets, there are infinitely many periodic orbits. 
By Lemma~\ref{lem:006}, we have that (4) holds.
\end{document}